\documentclass[acmtoms]{acmtrans2m_calc}

\usepackage{algorithmic}
\usepackage{algorithm}
\usepackage{psfrag}
\usepackage{epsfig}
\usepackage{amsmath}
\usepackage{color}
\usepackage{rotating}
\usepackage{pifont}
\usepackage{soul}
\usepackage{xspace}
\usepackage{subfigure}

\newcommand{\oh}[1]
    {\mbox{$ {\mathcal O}( #1 ) $}}
\newcommand{\ie}
    {{\em i.e.}~}
\newcommand{\eg}
    {{\em e.g.}~}
\newcommand{\ea}
    {{\em et al.}~}

\newcommand{\Checkmark}
    {\ding{51}\xspace}
\newcommand{\XSolid}
    {\ding{53}\xspace}
\newcommand{\dx}
    {\,\mbox{d}x}
\newcommand{\fig}[1]
    {Figure~\ref{fig:#1}\xspace}
\newcommand{\eqn}[1]
    {Equation~(\ref{eqn:#1})\xspace}
\newcommand{\sect}[1]
    {Section~\ref{sec:#1}\xspace}

\title{Increasing the Reliability of Adaptive Quadrature Using Explicit Interpolants}
            
\author{PEDRO GONNET\\ETH Z\"urich, Switzerland}

\begin{abstract}
    We present two new adaptive quadrature routines.
    Both routines differ from previously published algorithms
    in many aspects, most significantly in how they represent the integrand,
    how they treat non-numerical values of the integrand, how they deal with
    improper divergent integrals and how they estimate the integration error.
    The main focus of these improvements is to increase the {\em reliability}
    of the algorithms without significantly impacting their {\em efficiency}.
    Both algorithms are implemented in Matlab and tested using both the
    ``families'' suggested by Lyness and Kaganove and the battery
    test used by Gander and Gautschi and Kahaner. 
    They are shown to be more reliable, albeit in some cases
    less efficient, than other commonly-used adaptive integrators.
\end{abstract}


            

\begin{document}

\begin{bottomstuff} 
Author's address: P. Gonnet ({\tt gonnetp@inf.ethz.ch}), 
Department of Computer Science, ETH Zentrum, 8092 Z\"urich, Switzerland.
\end{bottomstuff}

\maketitle

\section{Introduction}
\label{sec:introduction}

Since the publication of the first adaptive quadrature algorithms almost 50
years ago, much has been done and even more has been written on the subject\footnote{
A recent review by the author \cite{ref:Gonnet2009} limited to error estimation listed 21 distinct
published algorithms and references to more than 50 publications
directly related to adaptive quadrature.}:
In 1962 Kuncir \cite{ref:Kuncir1962} kicked-off the field\footnote{
Davis and Rabinowitz \cite{ref:Davis1984} reference, as the first adaptive quadrature
routines, the works of Villars \cite{ref:Villars1956}, Henriksson \cite{ref:Henriksson1961}
and Kuncir \cite{ref:Kuncir1962}. Henriksson's algorithm, which appeared
in the first issue of BIT, is an ALGOL-implementation of the algorithm
described by Villars, which is itself an extension of an algorithm
by Morrin \cite{ref:Morrin1955}. These algorithms, however, are more reminiscent
of ODE integrators, which is why we will not consider them
to be ``genuine'' adaptive quadrature routines.} with his adaptive
Simpson's rule integrator, which uses -- as the name suggests -- Simpson's rule
to approximate the integral, bisecting recursively until the difference
between the approximation in one interval and that of its two sub-intervals
is below the required tolerance.
This approach, as simple as it may seem, still lives on, with some minor
modifications, as the default integrator {\tt quad} in {\small MATLAB}
(added in \cite{ref:Mathworks2003}).
{\tt quad} is itself a modification of Gander and Gautschi's 
2001 {\tt adaptsim} \cite{ref:Gander2001} which is itself a modification
of Lyness' 1970 {\tt SQUANK} \cite{ref:Lyness1970}, which is itself
almost identical to Kuncir's original algorithm\footnote{Lyness himself
formulates his algorithm as an extension to McKeeman's Adaptive Integrator
\cite{ref:McKeeman1962},
yet the resulting algorithm is much more similar to Kuncir's,
which was published merely a few months before McKeeman's.}.

At the time of Kuncir's publication,
McKeeman \cite{ref:McKeeman1962,ref:McKeeman1963,ref:McKeeman1963b}
and later Forsythe \ea \cite{ref:Forsythe1977}
extended this approach to use higher-degree Newton-Cotes rules
and/or sub-division into more than two sub-intervals.

In 1971, de~Boor \cite{ref:deBoor1971} 
introduced the concept of ``double adaptivity'', constructing
a Romberg T-table \cite{ref:Bauer1963} within each sub-interval to
approximate the integrand and deciding, at every step, whether to
extend the T-table by another row (\ie increase the order of the quadrature)
or to subdivide the interval.
This decision was made by using the convergence rates of the columns of
the T-table to guess the integrand's behavior, \ie ``well behaved'', singular, 
discontinuous or noisy, and apply an adequate strategy for that
behavior.

The emergence of more powerful computers and better algorithms 
(\eg \citeN{ref:Golub1969} and \citeN{ref:Gentleman1972b})
for the construction of more complex quadrature rules
quickly led to the wider use of Clenshaw-Curtis quadrature rules 
\cite{ref:Clenshaw1960}, as used by \citeN{ref:OHara1969}
and \citeN{ref:Oliver1972}, and later the use
of Gauss quadrature rules and their Kronrod extensions 
\cite{ref:Kronrod1965}, first used by \citeN{ref:Piessens1973}
and \citeN{ref:Patterson1973} independently.

Other interesting and/or noteworthy advances in the field are:
\begin{itemize}
    \item The introduction of stratified or recursively monotone stable (RMS) 
        quadrature rules \cite{ref:Sugiura1989,ref:Favati1991,ref:Laurie1992}, filling the gap
        between low-degree (due to low numerical stability at higher degrees)
        Newton-Cotes rules, the nodes of which
        nodes are re-usable over several recursion levels, and high-degree 
        (due to better numerical stability) yet non-reusable
        Clenshaw-Curtis or Gauss rules, thus providing
        some extra efficiency,
        
    \item The use of non-linear extrapolation when computing the integral
        or the error estimate \cite{ref:Rowland1972,ref:Venter2002,ref:Laurie1983,ref:deDoncker1978}, as is done
        in the highly successful {\tt QAGS} subroutine in the
        {\small QUADPACK} integration library,
        
    \item The use of higher-order coefficients relative to some base to
        compute the error estimate \cite{ref:OHara1968,ref:Oliver1972,ref:Berntsen1991}.
\end{itemize}

A number of authors have published comparisons of these and many other
adaptive quadrature routines \cite{ref:Casaletto1969,ref:Hillstrom1970,ref:Kahaner1971,ref:Malcolm1975,ref:Robinson1979,ref:Krommer1998,ref:Gonnet2009}
as well as methodologies to compare different routines \cite{ref:Lyness1977}.

As already noted in \citeN{ref:Rice1975},
despite all their differences, most adaptive quadrature algorithms follow
the general scheme, as in Algorithm~\ref{alg:general}.
First, an estimate of the integral in the interval $[a,b]$
is computed (Line~\ref{alg:general_q}).
An error estimate of the integral is then computed (in this example,
an absolute error estimate is approximated, Line~\ref{alg:general_eps}).
If this estimate is smaller than the required tolerance (Line~\ref{alg:general_if}),
then the estimate is returned (Line~\ref{alg:general_ret}).
Otherwise, the interval is bisected and the algorithm is called on both
halves (Line~\ref{alg:general_subdiv}) using a modified local tolerance $\tau'$.

\begin{algorithm}
    \caption{int $(f,a,b,\tau)$}
    \label{alg:general}
    \begin{algorithmic}[1]
        \STATE $Q \approx \int_a^b f(x)\,\mbox{d}x$ 
            \hfill ({\em approximate the integral in $[a,b]$}) \label{alg:general_q}
        \STATE $\varepsilon \approx \left| Q - \int_a^b f(x)\,\mbox{d}x \right|$
            \hfill ({\em approximate the integration error}) \label{alg:general_eps}
        \IF{$\varepsilon < \tau$} \label{alg:general_if}
            \RETURN{$Q$} \label{alg:general_ret}
                \hfill ({\em return the current estimate})
        \ELSE
            \RETURN{$\mbox{int}(f,a,(a+b)/2,\tau') + \mbox{int}(f,(a+b)/2,b,\tau')$}
                \\ \hfill ({\em call the integrator recursively on both sub-intervals}) \label{alg:general_subdiv}
        \ENDIF
    \end{algorithmic}
\end{algorithm}

\begin{algorithm}
    \caption{int $(f,a,b,\tau)$}
    \label{alg:general_nonrec}
    \begin{algorithmic}[1]
        \STATE $Q_0 \approx \int_a^b f(x)\,\mbox{d}x$ 
            \hfill ({\em approximate the integral in $[a,b]$})
        \STATE $\varepsilon_0 \approx \left| Q - \int_a^b f(x)\,\mbox{d}x \right|$
            \hfill ({\em approximate the integration error})
        \STATE $H \leftarrow \{[a,b,Q_0,\varepsilon_0]\}$ 
            \hfill ({\em initialize the heap with the first interval})
        \WHILE{$\sum_{\varepsilon_i \in H} \varepsilon_i > \tau$}
            \STATE $k \leftarrow \arg\max_k \varepsilon_k$
            \STATE $H \leftarrow H \setminus \{[a_k,b_k,Q_k,\varepsilon_k]\}$
                \hfill ({\em pop the interval with the largest error})
            \STATE $m \leftarrow (a_k+b_k)/2$
            \STATE $Q_\mathsf{left} \approx \int_{a_k}^m f(x)\,\mbox{d}x$
                \hfill ({\em compute the integral on the left})
            \STATE $\varepsilon_\mathsf{left} \approx \left| Q_\mathsf{left} - \int_{a_k}^m f(x)\,\mbox{d}x \right|$
                \hfill ({\em compute the error on the left})
            \STATE $Q_\mathsf{right} \approx \int_{m}^{b_k} f(x)\,\mbox{d}x$
                \hfill ({\em compute the integral on the right})
            \STATE $\varepsilon_\mathsf{right} \approx \left| Q_\mathsf{right} - \int_{m}^{b_k} f(x)\,\mbox{d}x \right|$
                \hfill ({\em compute the error on the right})
            \STATE $H \leftarrow H \cup \{ [a_k,m,Q_\mathsf{left},\varepsilon_\mathsf{left}] , [m,b_k,Q_\mathsf{right},\varepsilon_\mathsf{right}]\}$ 
                \\ \hfill ({\em push the new intervals back on the heap})
        \ENDWHILE
        \RETURN $\sum_{Q_i\in H}Q_i$
            \hfill ({\em return the sum of the integrals in the intervals})
    \end{algorithmic}
\end{algorithm}

Not all adaptive quadrature algorithms are recursive (locally adaptive):
many algorithms, such
as those in {\small QUADPACK}, maintain a heap of intervals and bisect the
interval with the largest local error estimate and return the new
sub-intervals to the heap until the sum of the local errors is below
the required tolerance (Algorithm~\ref{alg:general_nonrec}, globally adaptive).
This approach, although more memory-intensive, has several advantages
over the recursive approach, such as better control over the error estimate
and the ability to restart or refine an initial approximation
\cite{ref:Malcolm1975,ref:Rice1975}.

However, despite all these advances in numerical quadrature in general
and adaptive quadrature specifically, the results of these methods must
often be treated with caution, as failures are common even for relatively
simple integrands.
In this paper we will present two new adaptive quadrature algorithms
which attempt to address this lack of reliability.
The algorithms follow
the general scheme in Algorithm~\ref{alg:general_nonrec}, yet with 
significant differences to previous methods regarding how the integrand
is represented (Section~\ref{sec:function}), how the integration error
is estimated (Section~\ref{sec:error}) and how singularities
(Section~\ref{sec:singularities}) and
divergent integrals (Section~\ref{sec:divergent}) are treated.
The algorithm itself is presented in Section~\ref{sec:algorithm} and
in Section~\ref{sec:validation} it is validated against other popular
algorithms.
These results are then discussed in Section~\ref{sec:discussion}.

\section{Function Representation}
\label{sec:function}

In most quadrature algorithms, the integrand is not represented
internally except through different approximations of its integral.
We denote such approximations as
\begin{equation*}
    Q^{(m)}_n[a,b] = \sum_{i=1}^m Q_n[a + (i-1)h,a+ih] \approx \int_a^b f(x)\,\mbox{d}x, \quad h = \frac{b-a}{m}
\end{equation*}
where $n$ is the degree\footnote{In the following, we will use the term ``degree''
to specify the algebraic degree of precision of a quadrature rule, which is the
highest degree for which all polynomials of that degree will
always be integrated exactly by the rule.} of
the quadrature rule and $m$ its multiplicity.

The quadrature rule $Q_n$ itself is computed as the weighted sum of
the integrand evaluated at a pre-determined set of nodes\footnote{
In the following we assume, for notational simplicity, that the
number of nodes is the degree of the rule plus one. Although most quadrature
rules, \eg interpolatory quadrature rules with an odd number of symmetric
nodes or Gauss quadratures and their Kronrod extensions, need less than
$n+1$ nodes for degree $n$, this is a general upper bound for interpolatory
quadrature rules.} 
$x_i \in [-1,1]$, $i= 0 \dots n$:
\begin{equation}
    \label{eqn:qn}
    Q_n[a,b] = (b-a)\sum_{i=0}^n w_i f\left(\frac{a+b}{2} - \frac{a-b}{2}x_i\right).
\end{equation}

The evaluation of one or more such quadrature rules is usually the only
information considered regarding the integrand.

Some authors \cite{ref:Gallaher1967,ref:Ninomiya1980} use additional
nodes to numerically approximate the higher derivative directly using
divided differences, thus supplying additional information on the
integrand $f(x)$.
In a similar vein, \citeN{ref:OHara1969}, \citeN{ref:Oliver1972} and 
\citeN{ref:Berntsen1991} compute some of the higher-order coefficients
of the function relative to some orthogonal base, thus further
characterizing the integrand.

In all of these cases, however, the characterization of the integrand
is not complete and in most cases only implicit.
In the following, we will attempt to better characterize the integrand.

Before doing so, we note that for every
interpolatory quadrature rule, we are in fact computing a interpolating
polynomial $g_n(x)$ of degree $n$ such that
\begin{equation*}
    g_n(x_i) = f(x_i), \quad i=0 \dots n
\end{equation*}
and evaluating the integral thereof
\begin{equation*}
    Q_n[a,b] = \int_a^b g_n(x)\,\mbox{d}x.
\end{equation*}

This equivalence is easily demonstrated, as is done in many textbooks in numerical
analysis (\cite{ref:Stiefel1961,ref:Rutishauser1976,ref:Gautschi1997,ref:Schwarz1997,ref:Ralston1978} to name a few)\footnote{
If we consider the Lagrange interpolation $g_n(x)$ of the integrand
and integrate it, we obtain
\begin{equation*}
    \int_a^bg_n(x)\,\mbox{d}x = \int_a^b \sum_{i=0}^n \ell_i(x) f(x_i) \,\mbox{d}x
    = \sum_{i=0}^n f(x_i) \int_a^b \ell_i(x) \,\mbox{d}x = \sum_{i=0}^n f(x_i) w_i
\end{equation*}
where the $\ell_i(x)$ are the Lagrange polynomials and the
$w_i$ are the weights of the resulting quadrature rule.
}.

Since any polynomial interpolation of degree $n$ over $n+1$ distinct points is 
uniquely determined, it doesn't matter how we choose to represent $g_n(x)$ -- its
integral will always be identical to the result of
the interpolatory quadrature rule $Q_n[a,b]$ over the same nodes.

In the following, we will represent $g_n(x)$ as a linear combination of
orthogonal polynomials:
\begin{equation}
    \label{eqn:func_rep}
    g_n(x) = \sum_{i=0}^n c_i p_i(x)
\end{equation}
where the $p_i(x)$, $i=0 \dots n$ are polynomials
of degree $i$ which are orthonormal with respect to some inner product
\begin{equation*}
    (p_j,p_k) = \left\{ \begin{array}{ll}
        0\quad & j \ne k, \\ 1 & j = k. \end{array}\right.
\end{equation*}
We will use the coefficients $\mathbf{c} = (c_0,c_1,\dots,c_n)^\mathsf{T}$
from \eqn{func_rep} as our representation of $g_n(x)$.

For notational simplicity, we will assume that the integrand has been
transformed from the interval $[a,b]$ to the interval $[-1,1]$.
The polynomial $g_n(x)$ interpolates the integrand $f(x)$ at the nodes
$x_i \in [-1,1]$:
\begin{equation*}
    g_n(x_i) = f(x_i), \quad i=0 \dots n.
\end{equation*}

Given the function values $\mathbf{f} = (f(x_0),f(x_1),\dots,f(x_n))^\mathsf{T}$
at the nodes $x_i$, $i=0 \dots n$,
we can compute the coefficients by solving the linear system of equations
\begin{equation}
    \label{eqn:function_linsys}
    \mathbf{P} \mathbf{c} = \mathbf{f}
\end{equation}
where the matrix $\mathbf{P}$ with $P_{ij}=p_j(x_i)$
on the left-hand side is a {\em Vandermonde-like}
matrix.
The naive solution using Gaussian elimination is somewhat costly and may
be unstable \cite{ref:Gautschi1975}.
However, several algorithms exist to solve this problem stably in \oh{n^2} operations
for orthogonal polynomials satisfying a three-term recurrence relation
\cite{ref:Bjorck1970,ref:Higham1988,ref:Higham1990,ref:Gonnet2008b}.

In the following, we will use the orthonormal Legendre polynomials,
which are orthogonal with respect to the inner product
\begin{equation}
    \label{eqn:inner_product}
    (p_j,p_k) = \int_{-1}^1 p_j(x)p_k(x) \dx.
\end{equation}

We will evaluate and interpolate the integrand at the Chebyshev
nodes
\begin{equation*}
    x_i = \cos \left( \frac{\pi i}{n} \right), \quad i=0 \dots n.
\end{equation*}
These nodes are chosen over the Gauss quadrature nodes or equidistant
nodes due to their stability \cite{ref:Trefethen2008}, because
the nodes can be re-used when increasing the degree of the rule \cite{ref:Oliver1972}
and because they include the interval boundaries.

The resulting Vandermonde-like matrix has a condition number
$\kappa_\infty(\mathbf{P}) \in \mathcal O(n^{3/2})$
which remains $<1\,000$ for $n \leq 100$ and is thus tractable
even for moderate $n$ \cite{ref:Gonnet2009}.

The resulting representation of $g_n(x)$ (\eqn{func_rep}) has some interesting properties.
First of all, it is simple to evaluate the integral of $g_n(x)$ using
\begin{equation}
    \label{eqn:coeffint}
    \int_{-1}^1 g_n(x) \,\mbox{d}x = \int_{-1}^1 \sum_{i=0}^n c_i p_i(x) \,\mbox{d}x
    = \sum_{i=0}^n c_i \underbrace{\int_{-1}^1 p_i(x) \,\mbox{d}x}_{=\omega_i}
    = \boldsymbol\omega^\mathsf{T} \mathbf{c}
\end{equation}
where the weights $\boldsymbol\omega^\mathsf{T}$ can be pre-computed
and applied much in the same way as the weights of a quadrature rule.
Note that for the normalized Legendre polynomials used herein,
$\mathbf{\omega}^\mathsf{T} = ( 1/\sqrt{2} , 0 , \dots , 0 )$.

We can also evaluate the $L_2$-norm of $g_n(x)$ quite efficiently using
Parseval's theorem
\begin{equation*}
    \left[ \int_{-1}^1 g_n^2(x) \,\mbox{d}x \right]^{1/2} =
        \left[ \sum_{i=0}^n c_i^2 \right]^{1/2} = \|\mathbf{c}\|_2.
\end{equation*}
In the following, we will use $\|\cdot\|$ to denote the 2-norm
for vectors.

A final useful feature is that, given the coefficients of $g_n(x)$
on $[-1,1]$, we can construct upper-triangular matrices
\begin{equation*}
    T^{(\ell)}_{i,j} = \int_{-1}^1 p_i(x)p_j\left(\frac{x-1}{2}\right)\,\mbox{d}x, \quad
    T^{(r)}_{i,j} = \int_{-1}^1 p_i(x)p_j\left(\frac{x+1}{2}\right)\,\mbox{d}x, \quad i = 0\dots n, j \ge i
\end{equation*}
such that
\begin{equation}
    \label{eqn:tlr}
    \mathbf{c}^{(\ell)} = \mathbf{T}^{(\ell)} \mathbf{c} \quad \mbox{and} \quad \mathbf{c}^{(r)} = \mathbf{T}^{(r)} \mathbf{c}
\end{equation}
are the coefficients of $g_n(x)$ on the left and right sub-intervals
$[-1,0]$ and $[0,1]$ respectively.
These matrices depend only on the polynomials $p_i(x)$ and can therefore
be pre-computed for any set of nodes such that\footnote{
This can be shown by representing the polynomials $p_i((x-1)/2)$ of degree $i$
as a linear combination of the polynomials $p_j(x)$, $j=0\dots i$ where
the coefficients are computed using the inner product in \eqn{inner_product}:
\begin{equation*}
    p_i\left( \frac{x-1}{2} \right) \quad = \quad \sum_{j=0}^i p_j(x)\left[ \int_{-1}^1 p_j(x)p_i\left(\frac{x-1}{2}\right)\,\mbox{d}x \right] \quad = \quad \sum_{j=0}^i p_j(x) T^{(\ell)}_{j,i}.
\end{equation*}
We can then re-insert this expression into \eqn{ggn}
\begin{equation*}
    g^{(\ell)}_n(x) = g_n\left( \frac{x-1}{2} \right) \quad \Longrightarrow \quad
    \sum_{i=0}^nc^{(\ell)}_ip_i(x) = \sum_{i=0}^nc_ip_i\left( \frac{x-1}{2} \right) = 
    \sum_{i=0}^n c_i \sum_{j=0}^i p_j(x) T^{(\ell)}_{j,i}
\end{equation*}
which, swapping the indices $i$ and $j$ and re-arranging the sums
on the right-hand side can be re-written as
\begin{equation*}
    \sum_{i=0}^nc^{(\ell)}_ip_i(x) = \sum_{i=0}^n p_i(x) \sum_{j=i}^n c_j T^{(\ell)}_{i,j}
\end{equation*}
using which the vector of coefficients $\mathbf c^{(\ell)}$ can be computed
as in \eqn{tlr}.
}
\begin{equation}
    \label{eqn:ggn}
    g^{(\ell)}_n(x) = g_n\left( \frac{x-1}{2} \right), \quad x \in [-1,1].
\end{equation}
The coefficients $\mathbf{c}^{(\ell)}$ and $\mathbf{c}^{(r)}$
can be useful if, after bisecting
an interval, we want to re-use, inside one of the sub-intervals, the
interpolation computed over the entire original interval.

\section{Error Estimation}
\label{sec:error}

This section contains a summary of the more important results
of \cite{ref:Gonnet2009}.
For a more complete discussion and testing of the error estimators
presented herein,
we refer to that publication.

\begin{figure}
    \centerline{\epsfig{file=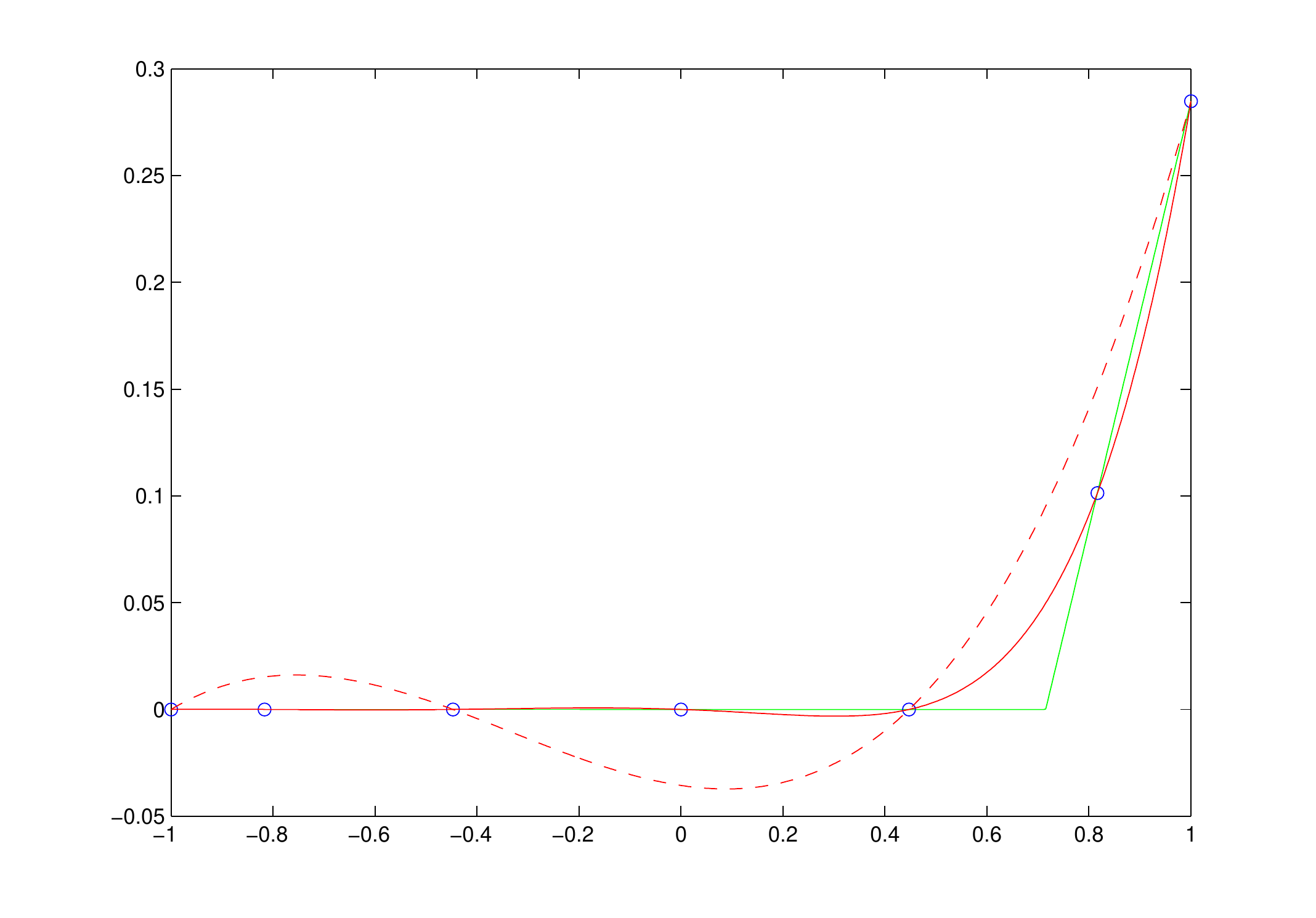,width=0.8\textwidth}}
    \caption{A discontinuous function (solid) and the integrand
        interpolation using the Gauss-Lobatto rule $Q_{5}[-1,1]$ (dashed)
        and the Gauss-Kronrod rule $Q_{9}[-1,1]$ (dotted)
        from Matlab's {\tt quadl} integrator.
        Note that although both quadratures
        return the same result, the interpolated functions differ
        significantly.}
    \label{fig:discont}
\end{figure}

Although they differ in their specific implementations,
what all these error estimates have in common is that they try to approximate
the quantity
\begin{equation}
    \label{eqn:err_old}
    \varepsilon = \left| Q^{(m)}_n[a,b] - \int_a^bf(x)\,\mbox{d}x \right|
\end{equation}
using only two or more approximations of the integral or of its
coefficients relative to some base.
In these estimates, problems may occur when the difference
between two estimates $Q^{(m_1)}_n[a,b]$ and $Q^{(m_2)}_n[a,b]$
or $Q_{n_1}[a,b]$ and $Q_{n_2}[a,b]$, or the magnitude
of the computed coefficients is {\em accidentally small}\,\footnote{
This term was first used by \citeN{ref:OHara1968} to describe this problem.},
\ie the approximations used to compute the error estimate
are too imprecise, resulting in a false small error estimate.
This is often the case near singularities and discontinuities where the
assumptions on which the error estimate is based, \eg continuity and/or
smoothness, do not hold.
If we re-construct the underlying interpolatory polynomials for the
pair of quadrature rules used, we see that the interpolations
differ significantly (\eg see \fig{discont}).
This difference is a good indicator for whether the integrand
is correctly represented or not.

It is for this reason that instead of trying to compute the error as
in \eqn{err_old}, we will use the $L_2$ norm\footnote{
Note that if instead of using the orthonormal Legendre polynomials
we were to use polynomials orthogonal with respect to any specific
measure $w(x)$, we would compute the $L_2$-norm
with respect to that measure:
\begin{equation*}
    \varepsilon = \left[ \int_{-1}^1 w(x)(g_n(x)-f(x))^2\dx \right]^{1/2}.
\end{equation*}
Fortunately enough, for any measure $w(x)$, the following derivations apply
without modification.
}
of the
difference between the integrand $f(x)$ and the interpolant $g_n(x)$
\begin{equation}
    \label{eqn:err_new}
    \left[ \int_{-1}^1 \left( g_n(x) - f(x) \right)^2 \,\mbox{d}x \right]^{1/2}.
\end{equation}
This is an approximation
of the integration error \eqn{err_old}.
The error estimate will only be zero if the interpolated
integrand matches the integrand on the entire interval
\begin{equation*}
    g_n(x) = f(x), \quad x \in [-1,1].
\end{equation*}
In such a case, the integral will also be computed exactly.
The error \eqn{err_new} is therefore, assuming we can evaluate
it reliably, not susceptible to ``accidentally small'' values.

Since we do not know $f(x)$ explicitly, \ie we can only sample
$f(x)$ in a point-wise fashion, we cannot evaluate the
right-hand side of \eqn{err_new} exactly.
In a first, naive approach, we could compute two interpolations
$g^{(1)}_{n_1}(x)$ and $g^{(2)}_{n_2}(x)$ of different degree
where $n_1 < n_2$.
If we assume, as is done for error estimators using quadrature
rules of differing degree, that $g^{(2)}_{n_2}(x)$ is a sufficiently
precise approximation of the integrand
\begin{equation*}
    g^{(2)}_{n_2}(x) \approx f(x), \quad x \in [-1,1]
\end{equation*}
then we can approximate the error of the interpolation $g^{(1)}_{n_1}(x)$
as
\begin{eqnarray*}
    \left[ \int_{-1}^1 \left( f(x) - g^{(1)}_{n_1}(x) \right)^2 \,\mbox{d}x \right]^{1/2}
        & \approx & \left[ \int_{-1}^1 \left( g^{(2)}_{n_2}(x) - g^{(1)}_{n_1}(x) \right)^2 \,\mbox{d}x \right]^{1/2} \\
        & = & \| \mathbf{c}^{(2)} - \mathbf{c}^{(1)} \|
\end{eqnarray*}
where $\mathbf{c}^{(1)}$ and $\mathbf{c}^{(2)}$ are the vectors
of the coefficients of $g^{(1)}_{n_1}(x)$ and $g^{(2)}_{n_2}(x)$
respectively and $c^{(1)}_i = 0$ where $i > n_1$.
Our first, naive error estimate is hence
\begin{equation}
    \label{eqn:err_naive}
    \varepsilon_\mathsf{naive} := \|\mathbf{c}^{(1)}-\mathbf{c}^{(2)}\|.
\end{equation}

This error estimate, however, only applies to the lower-degree 
estimate $g^{(1)}_{n_1}(x)$.
Yet if we are going to compute the higher-degree estimate $g^{(2)}_{n_2}(x)$,
it would be preferable to have an error estimate for {\em that} approximation.

We can, taking a different approach, use the interpolation error
\begin{equation}
    \label{eqn:err_interp}
    |g_n(x) - f(x)| = \left|\frac{f^{(n+1)}(\xi_x)}{(n+1)!}\pi_n(x)\right| , \quad x \in [-1,1]
\end{equation}
where $\xi_x \in [-1,1]$ depends on the value of $x$ and where
$\pi_n(x) = \prod_{i=0}^n (x - x_i)$
is the Newton basis polynomial over the nodes of the interpolation $g_n(x)$.
Taking the $L_2$-norm on both sides of \eqn{err_interp}
we obtain
\begin{equation*}
    \varepsilon = \left[ \int_{-1}^1 \left( g_n(x) - f(x) \right)^2 \,\mbox{d}x \right]^{1/2} = \left[ \int_{-1}^{1} \left(\frac{f^{(n+1)}(\xi_x)}{(n+1)!}\right)^2 \pi^2_n(x) \,\mbox{d}x \right]^{1/2}.
\end{equation*}
Since $\pi_n^2(x)$ is, by definition, positive for any $x$,
we can apply the mean value theorem of integration and extract
the derivative resulting in
\begin{equation}
    \label{eqn:err_final}
    \varepsilon = \left[ \int_{-1}^1 \left( g_n(x) - f(x) \right)^2 \,\mbox{d}x \right]^{1/2} = \left|\frac{f^{(n+1)}(\xi)}{(n+1)!}\right|\left[ \int_{-1}^{1} \pi^2_n(x) \,\mbox{d}x \right]^{1/2}, \quad \xi \in [-1,1].
\end{equation}

Given two interpolations $g^{(1)}_n(x)$ and $g^{(2)}_n(x)$ over a
non-identical set of nodes, we can compute the interpolation errors
\begin{equation}
    \left|g^{(\star)}_n(x) - f(x)\right| = \left|\frac{f^{(n+1)}(\xi_\star)}{(n+1)!} \pi^{(\star)}_n(x)\right|, \quad \xi_\star \in [-1,1], \quad \star \in \{1,2\} \label{eqn:err_g1g2}
\end{equation}
where $\pi^{(1)}_n(x)$ and $\pi^{(2)}_n(x)$ are the Newton basis polynomials
over the nodes of $g^{(1)}_n(x)$ and $g^{(2)}_n(x)$ respectively.
If we assume that $f^{(n+1)}(x)$ is constant for $x\in[-1,1]$
(a stricter version of the ``sufficiently smooth'' assumption for the
purpose of deriving this error estimate)
and take the $L_2$-norm of the difference between both errors, we obtain
\begin{equation}
    \label{eqn:err_new2}
    \left[ \int_{-1}^1 \left( g^{(1)}_n(x) - g^{(2)}_n(x) \right)^2\,\mbox{d}x \right]^{1/2} =
    \left|\frac{f^{(n+1)}(\xi)}{(n+1)!}\right| \left[ \int_{-1}^1 \left( \pi^{(1)}_n(x) - \pi^{(2)}_n(x) \right)^2\,\mbox{d}x \right]^{1/2}.
\end{equation}

If we represent the interpolations $g^{(1)}_n(x)$ and $g^{(2)}_n(x)$
by their coefficients $\mathbf{c}^{(1)}$ and $\mathbf{c}^{(2)}$
respectively, then we can write the left-hand side of
\eqn{err_new2} as
\begin{equation}
    \label{eqn:err_new3}
    \left\| \mathbf{c}^{(1)} - \mathbf{c}^{(2)} \right\| =
    \left|\frac{f^{(n+1)}(\xi)}{(n+1)!}\right| \left[ \int_{-1}^1 \left( \pi^{(1)}_n(x) - \pi^{(2)}_n(x) \right)^2\,\mbox{d}x \right]^{1/2}.
\end{equation}
Similarly, if we represent the Newton basis polynomials $\pi^{(1)}_n(x)$ 
and $\pi^{(2)}_n(x)$ by their coefficients $\mathbf{b}^{(1)}$ and
$\mathbf{b}^{(2)}$ respectively
\begin{equation}
    \label{eqn:err_newton}
    \pi^{(1)}_n(x) = \sum_{i=0}^{n+1} b^{(1)}_i p_i(x), \quad \pi^{(2)}_n(x) = \sum_{i=0}^{n+1} b^{(2)}_i p_i(x),
\end{equation}
we can isolate the fraction on
the right hand side of \eqn{err_new3}
\begin{equation}
    \label{eqn:err_new4}
    \left|\frac{f^{(n+1)}(\xi)}{(n+1)!}\right| = 
        \frac{\left\| \mathbf{c}^{(1)} - \mathbf{c}^{(2)} \right\|}{\left\| \mathbf{b}^{(1)} - \mathbf{b}^{(2)} \right\|}.
\end{equation}

Inserting this expression into the original error estimate
(\eqn{err_final}) for the interpolation $g^{(1)}_n(x)$ 
we then obtain
\begin{equation*}
    \left[ \int_{-1}^1 \left( g^{(1)}_n(x) - f(x) \right)^2\,\mbox{d}x \right]^{1/2} 
        = \displaystyle \|\mathbf{b}^{(1)}\|\frac{\left\| \mathbf{c}^{(1)} - \mathbf{c}^{(2)} \right\|}{\left\| \mathbf{b}^{(1)} - \mathbf{b}^{(2)} \right\|}
        =: \varepsilon^{(1)} \\
\end{equation*}
Hence, using two interpolations of the same degree, we obtain
the more refined error estimate
\begin{equation}
    \label{eqn:err_refined}
    \varepsilon_\mathsf{ref} := \frac{\|\mathbf{c}^{(1)}-\mathbf{c}^{(2)}\|}{\|\mathbf{b}^{(1)}-\mathbf{b}^{(2)}\|}\|\mathbf{b}^{(1)}\|
\end{equation}
for the interpolation $g^{(1)}_n(x)$.

Note that if the nodes of the interpolations $g^{(1)}_n(x)$ 
and $g^{(2)}_n(x)$ are fixed, we can pre-compute the scaling
$\|\mathbf{b}^{(1)}\|/\|\mathbf{b}^{(1)}-\mathbf{b}^{(2)}\|$.

Instead of explicitly computing two different interpolations over an interval
$[a,b]$ to construct the error estimate, we can re-use the interpolation
from the previous level of recursion after bisection, the coefficients
$\mathbf{c}^\mathsf{old}$ of which can be computed using \eqn{tlr}.
Likewise, we can compute the coefficients $\mathbf{b}^\mathsf{old}$
of the Newton basis polynomial over the nodes of the previous level in the same way.
However, since $\mathbf{b}^\mathsf{old}$ and $\mathbf{b}$ are not in the same
interval, we have to scale the coefficients of $\mathbf{b}^\mathsf{old}$
by $2^{n+1}$ such that \eqn{err_g1g2} holds.
In this way, we only need to compute and store a single matrix
$\mathbf{P}^{-1}$ and vector $\mathbf{b}$ for a single stencil of
interpolation nodes.

As with the previous error estimators, we have also made an assumption 
of smoothness regarding the integrand by assuming that
$f^{(n+1)}(x)$ is constant for $x\in[-1,1]$
to construct \eqn{err_refined}.
We can't verify this directly, but we can verify if our computed
$|\frac{f^{(n+1)}(\xi)}{(n+1)!}|$ (\eqn{err_new4})
actually satisfies \eqn{err_g1g2} for the nodes of the first
interpolation by testing if
\begin{equation}
    \label{eqn:err_test}
    \left| g^{(2)}_n(x_i) - f(x_i) \right| \leq \vartheta_1 \left|\frac{f^{(n+1)}(\xi)}{(n+1)!}\right| \left|\pi^{(2)}_n(x_i)\right|
\end{equation}
is satisfied for all $i=0 \dots n$, where the $x_i$ are the nodes of
the interpolation $g^{(1)}_n(x)$.
The value $\vartheta_1 \geq 1$ is an arbitrary relaxation parameter
(for the tests in \sect{validation} we use $\vartheta_1=1.1$).
If this condition is violated for any of the $x_i$, then we use
the naive error estimate in \eqn{err_naive}.

\section{Singularities and Undefined Values}
\label{sec:singularities}

Since most adaptive quadrature algorithms are designed for general-purpose
use, they will often be confronted with
integrands containing singularities or undefined function values.
These can cause problems on two levels:
\begin{itemize}
    \item The quadrature rule has to be evaluated with a non-numerical
        value such as a $\mathsf{NaN}$ or $\pm\mathsf{Inf}$,
    \item The integrand may not be as smooth and continuous
        as the algorithm might assume.
\end{itemize}

Such problems arise when integrating functions such as
\begin{equation*}
    \int_0^h x^\alpha \dx, \quad \alpha < 0
\end{equation*}
which have a singularity at $x=0$, or when computing seemingly
innocuous integrals such as
\begin{equation*}
    \int_0^h \frac{\sin x}{x} \dx
\end{equation*}
for which the integrand is undefined at $x=0$, yet has a
well-defined limit
\begin{equation*}
    \lim_{x\rightarrow 0} \frac{\sin x}{x} = 1.
\end{equation*}

In both cases problems {\em could} be avoided by either shifting
the integration domain slightly or by modifying the integrand such
as to catch the undefined cases and return a correct numerical result.
This would, however, require some prior reflection and intervention by
the user, which would defeat the purpose of a general-purpose 
quadrature algorithm.

Most algorithms deal with singularities by ignoring them, setting the offending
value of the integrand to 0 \cite[Section 2.12.7]{ref:Davis1984}.
Another approach, taken by {\tt quad} and {\tt quadl} in Matlab,
is to shift the edges of the domain by $\varepsilon_\mathsf{mach}$ if
a non-numerical value is encountered there and to abort with a warning
if a non-numerical value is encountered elsewhere in the interval.
Since singularities may exist explicitly at the boundaries
(\eg integration of $x^\alpha$, $\alpha < 0$ in the range $[0,h]$),
the explicit treatment of the 
boundaries is needed, whereas for arbitrary singularities within the interval,
the probability of hitting them exactly is somewhat small.

{\small QUADPACK}'s {\tt QAG} and {\tt QAGS} algorithms take
a similar approach: since the nodes of the Gauss and Gauss-Lobatto
quadrature rules used therein do not include the interval 
boundaries, non-numerical values at the interval boundaries will
be implicitly avoided.
If the algorithm has the misfortune of encountering such a value
inside the interval, it aborts.

Our approach to treating singularities will
be somewhat different: instead of setting non-numerical values
to 0, we will simply {\em remove} that node from our interpolation
of the integrand.
This can be done rather efficiently by computing the interpolation
as shown before using a function value of $f(x_j)=0$ for the offending
$j$th node and then {\em down-dating} (as opposed to up-dating)
the interpolation, \ie
removing the $j$th node from the interpolation, resulting in an 
interpolation of degree $n-1$:
\begin{equation*}
    g_{n-1}(x) = \sum_{i=0}^{n-1} c^{(n-1)}_i p_i(x)
\end{equation*}
which still interpolates the integrand at the remaining $n$ nodes.

The coefficients $c^{(n-1)}_i$ of $g_{n-1}(x)$ can be computed,
as described in \cite{ref:Gonnet2008b}, using
\begin{equation*}
    \label{eqn:sing_downdatec}
    c^{(n-1)}_i = c_i - \frac{c_n}{b^{(n-1)}_n}b^{(n-1)}_i, \quad i = 0 \dots n
\end{equation*}
where the $c_i$ are the computed coefficients of $g_n(x)$
and the $b^{(n-1)}_i$ are the coefficients of the downdated
Newton polynomial computed by solving the upper-triangular
system of equations
\begin{equation*}
    \label{eqn:sing_downdateb}
    \left( \begin{array}{ccccc}
        \alpha_0 & -(x_j + \beta_1) & \gamma_2 \\
        & \ddots & \ddots & \ddots \\
        & & \alpha_{n-2} & -(x_j + \beta_{n-1}) & \gamma_n \\
        & & & \alpha_{n-1} & -(x_j + \beta_n) \\
        & & & & \alpha_n
    \end{array}\right)
    \left( \begin{array}{c} b^{(n-1)}_0 \\ b^{(n-1)}_1 \\ \vdots \\ b^{(n-1)}_{n-1} \end{array}\right)
    =
    \left( \begin{array}{c} b_1 \\ b_2 \\ \vdots \\ b_n \end{array}\right)
\end{equation*}
using back-substitution, where the $b_i$ are the coefficients
of the Newton polynomial over the nodes of the quadrature
rule (\eqn{err_newton}) and the $\alpha_i$, $\beta_i$ and
$\gamma_i$ are the coefficients of the three-term recurrence
relation satisfied by the polynomials of the orthogonal basis:
\begin{equation*}
    \alpha_k p_{k+1}(x) = (x + \beta_k)p_k(x) - \gamma_k p_{k-1}(x).
\end{equation*}

The modified vectors $\mathbf{c}^{(n-1)}$ and
$\mathbf{b}^{(n-1)}$ are then used in the same way as
$\mathbf{c}$ and $\mathbf{b}$ respectively for the 
computation of the integral and of the error estimate.

\section{Divergent Integrals}
\label{sec:divergent}

Divergent integrals are integrals which tend to
$\pm \infty$ and thus cause most algorithms to either recurse
infinitely or return an incorrect finite result.
They are usually caught by limiting the
recursion depth or the number of function
evaluations artificially.
Both approaches do not {\em per se} attempt to detect divergent
behavior, and may therefore cause the algorithm to fail for
complicated yet non-divergent integrals.

\citeN{ref:Ninham1966} studied the approximation error when computing
\begin{equation}
    \label{eqn:divergent}
    \int_0^h x^\alpha\,\mbox{d}x
\end{equation}
using the trapezoidal rule.
The integral exists for $\alpha > -1$ and is divergent otherwise.
Following his analysis, we compute the refined error estimate described
in \sect{error} (\eqn{err_refined})
for the intervals $[0,h]$ and $[0,h/2]$
using an 11-node Clenshaw-Curtis quadrature rule\footnote{
Remember that we are interpolating over the Chebyshev nodes which
is equivalent to using a Clenshaw-Curtis quadrature rule.} and removing the
singular node at $x=0$ as described above.

\begin{figure}
    \centerline{\epsfig{file=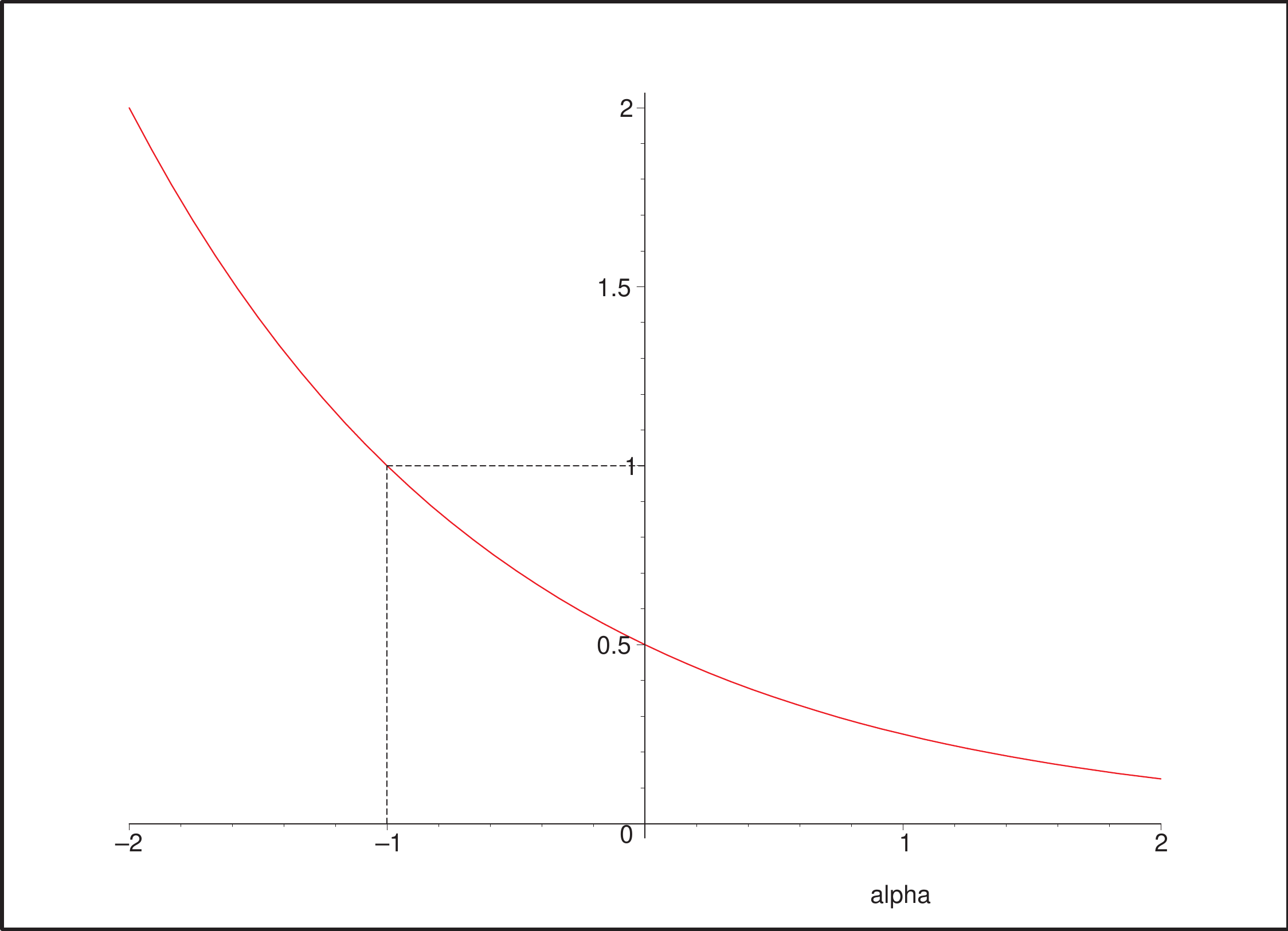,width=0.5\textwidth,angle=-90}}
    \caption{Ratio of the error estimates for $\int_0^hx^\alpha\,\mbox{d}x$
        for the intervals $[0,h/2]$ over $[0,h]$ for different $\alpha$.
        Note that the error grows (ratio $>1$) for $\alpha < -1$,
        where the integral is divergent.}
    \label{fig:alpha}
\end{figure}

We note that as the algorithm recurses to the leftmost interval,
the local error estimate as well as the computed integral itself
remain constant for $\alpha = -1$ and 
{\em increase} for $\alpha < -1$.
In Figure~\ref{fig:alpha} we plot the ratio of the error estimate in
the left sub-interval over the error in the entire interval,
$\varepsilon[0,h/2]/\varepsilon[0,h]$, over the parameter $\alpha$.
For $\alpha = -1$ the ratio is $1$ meaning that the error estimate of
the leftmost interval remains constant even after halving the
interval.
For $\alpha<-1$, for which the integral diverges, the error estimate
in the left half-interval is {\em larger} than the error
estimate over the entire interval.

\begin{figure}
    \centerline{\epsfig{file=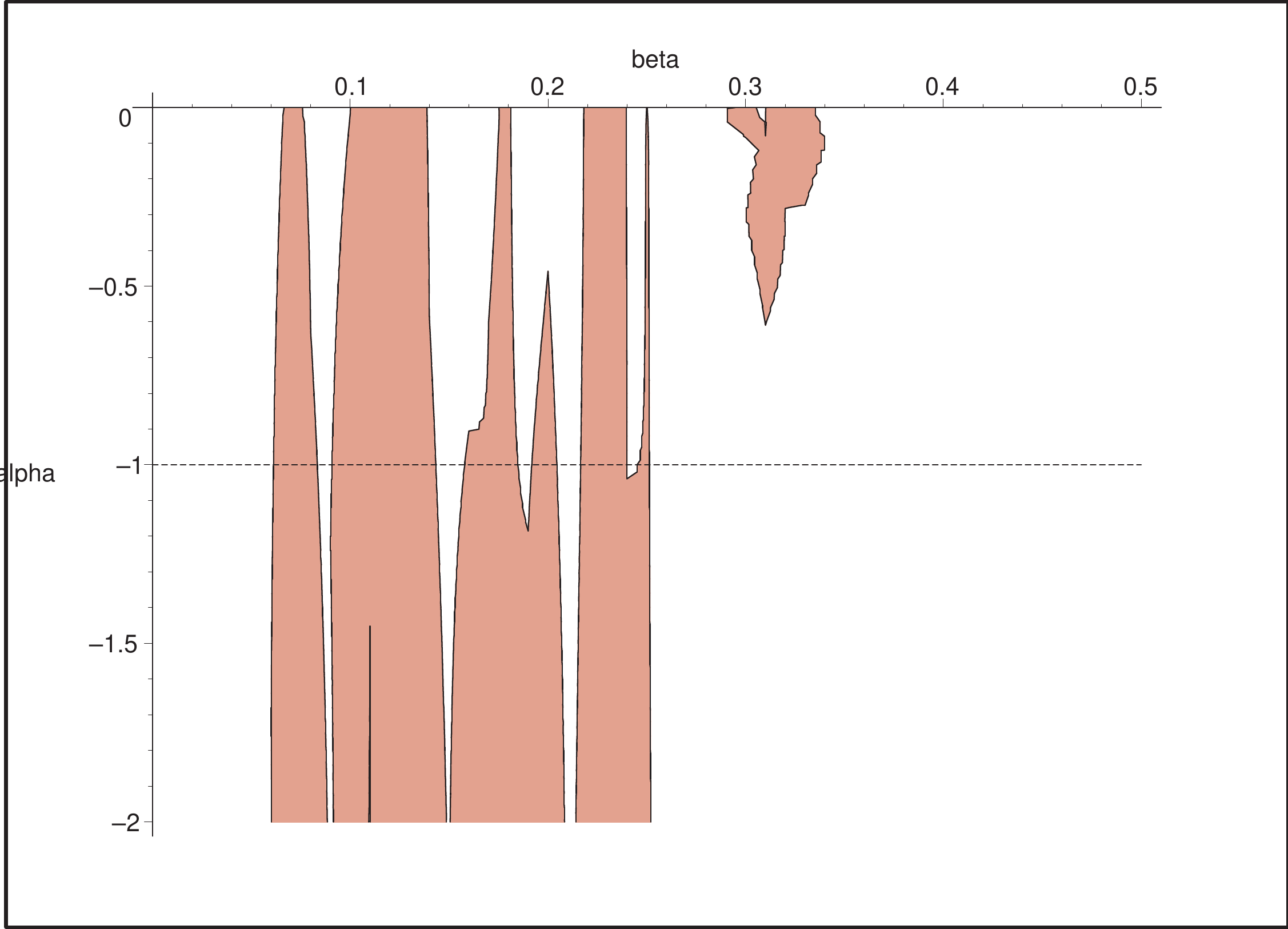,width=0.5\textwidth,angle=-90}}
    \caption{Contour of the ratio of the error estimates for
        $\int_0^h|x-\beta|^\alpha\,\mbox{d}x$ over the intervals $[0,1]$
        and $[0,1/2]$. The filled area represent the region in
        which this ratio is larger than 1.}
    \label{fig:alphabeta}
\end{figure}

The rate at which the error decreases (or, in this case, increases)
may be a good indicator for
the convergence or divergence of the integral in \eqn{divergent}, where
the singularity is at the edge of the domain, yet
it does not work as well for the shifted singularity
\begin{equation}
    \label{eqn:divergent2}
    \int_0^h |x-\beta|^\alpha\,\mbox{d}x, \quad \beta \in [0,h/2].
\end{equation}
Depending on the location of the singularity ($x=\beta$),
the ratio of the error estimates over $[0,h]$ and $[0,h/2]$
varies widely for both $\alpha > -1$ and $\alpha \leq -1$
and can not be used to determine whether the
integral diverges or not.
In \fig{alphabeta} we have shaded the regions in which the ratio
of the error estimates $\varepsilon[0,h/2]/\varepsilon[0,h] > 1$ for
different values of $\alpha$ and the location of the singularity $\beta$.
For this ratio to be a good indicator for the value of $\alpha$
(and hence the convergence/divergence of the integral), the shaded area
should at least partially cover the lower half of the plot where
$\alpha < -1$, which it does not.

\begin{figure}
    \centerline{\epsfig{file=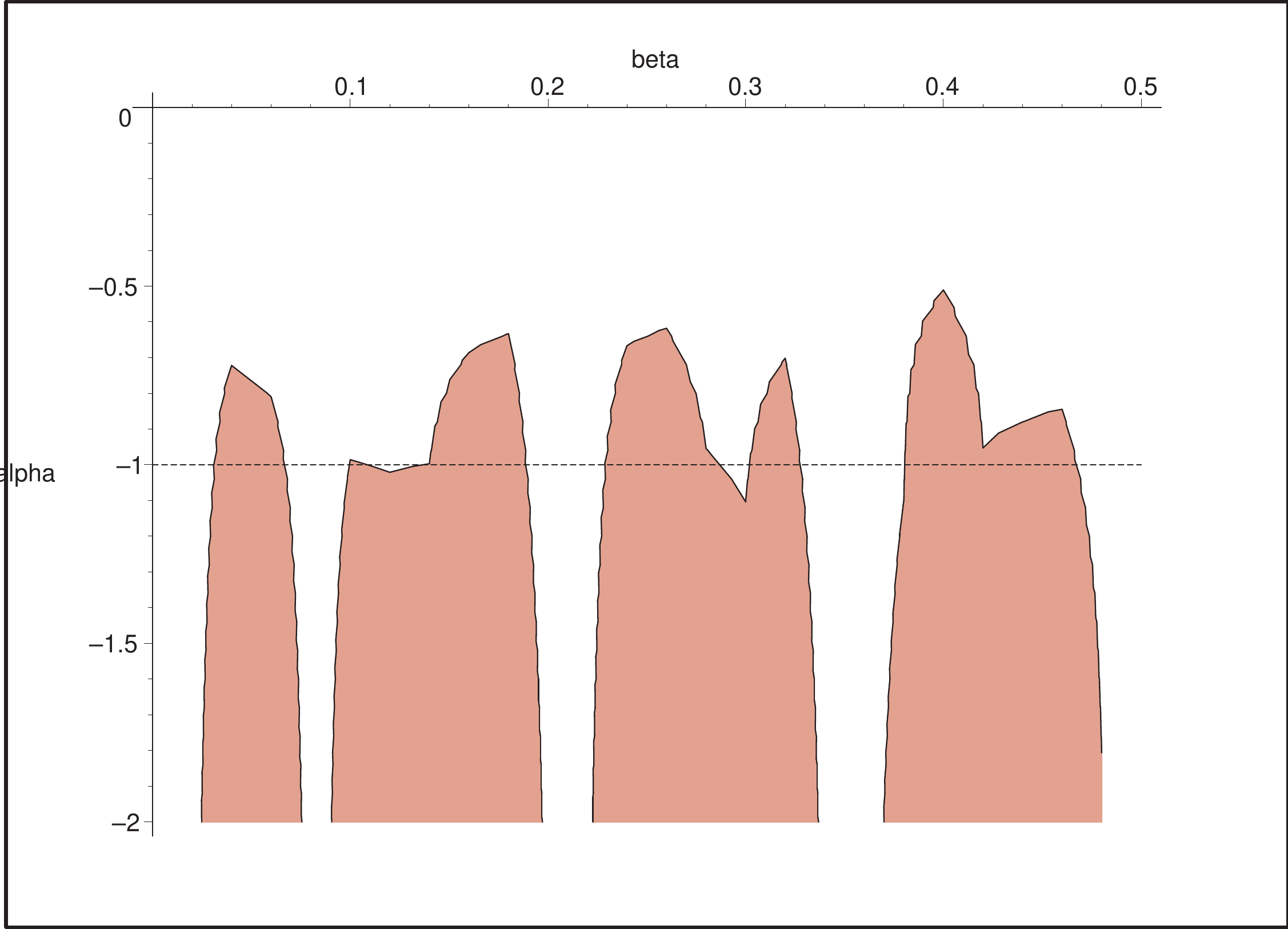,width=0.5\textwidth,angle=-90}}
    \caption{Contour of the ratio of the integral estimates for
        $\int_0^h|x-\beta|^\alpha\,\mbox{d}x$ over the intervals $[0,1]$
        and $[0,1/2]$. The filled area represent the region in
        which this ratio is larger than 1.}
    \label{fig:alphabeta2}
\end{figure}

A more reliable approach consist of comparing the computed
{\em integral} in two successive intervals $[a,b]$ and $[a,(a+b)/2]$
or $[(a+b)/2,b]$.
For the integrand in \eqn{divergent}, the integral
in the left sub-interval $[0,h/2]$ is {\em larger} than that over
the interval $[0,h]$ when $\alpha < -1$.
For the integral in \eqn{divergent2} the ratio
of the integrals in the intervals $[0,h]$ and $[0,h/2]$ is larger
than 1 for {\em most} cases where $\alpha \leq -1$ (see
Figure~\ref{fig:alphabeta2}).

Although this relation (ratio $> 1 \Rightarrow \alpha < -1$),
which is independent of the interval $h$, is not always correct, 
it can still be used as a statistical {\em hint} for the integrand's behavior
during subdivision.
In the area $-2 \leq \alpha \leq -1$ it is correct approximately
two thirds of the time.
We will therefore count the number of times that
\begin{equation}
    \label{eqn:sing_test}
    \frac{Q_n[a,(a+b)/2]}{Q_n[a,b]} \geq 1 \quad \mbox{or} \quad
    \frac{Q_n[(a+b)/2,b]}{Q_n[a,b]} \geq 1
\end{equation}
during subdivision.
Note that since the integral goes to either $+\infty$ or $-\infty$,
the integral approximations will be of the same sign and thus the sign
of the ratios do not matter.
If this count exceeds some maximum number and is more than half of the
recursion depth -- \ie the ratio of integrals was larger than one over
more than half of the subdivisions -- then we
declare the integral to be divergent and return an error or warning
to the user to this effect.

\section{The Algorithm}
\label{sec:algorithm}

We present the algorithm in two variants, Algorithm~\ref{alg:final_naive}
and Algorithm~\ref{alg:final_refined}, using both the naive and the
refined error estimates presented in Section~\ref{sec:error}.
Both algorithms follow the globally adaptive general scheme shown in 
Algorithm~\ref{alg:general_nonrec}.

\begin{algorithm}
    \caption{int\_naive $(f,a,b,\tau)$}
    \label{alg:final_naive}
    \begin{algorithmic}[1]
    
        \IFOR{$i=0 \dots n_{d_\mathsf{max}}$}{$f_i \leftarrow f\left((a+b)/2 - (a-b)x^{(m)}_i/2\right)$}
            \\ \hfill ({\em evaluate the integrand at the nodes $x^{(m)}_i$})\label{alg:naive_init1}
        
        \STATE $\mathbf{c}^{(d_\mathsf{max}-1)} \leftarrow (\mathbf{P}^{(d_\mathsf{max}-1)})^{-1} \mathbf{f}[1:2:n_m+1]$,
            $\mathbf{c}^{(d_\mathsf{max})} \leftarrow (\mathbf{P}^{(d_\mathsf{max})})^{-1} \mathbf{f}$
            \\ \hfill ({\em compute the interpolation coefficients})
            
        \STATE $q_0 \leftarrow (b - a) \boldsymbol\omega^\mathsf{T}\mathbf{c}^{(d_\mathsf{max})}/2$
            \hfill ({\em approximate the integral as per \eqn{coeffint}})
        
        \STATE $\varepsilon_0 \leftarrow (b - a) \|\mathbf{c}^{(d_\mathsf{max})}-\mathbf{c}^{(d_\mathsf{max}-1)}\|/2$
            \hfill ({\em approximate the error})\label{alg:naive_init2}
        
        \STATE $H \leftarrow \{[a,b,\mathbf{c}^{(d_\mathsf{max})},q_0,\varepsilon_0,d_\mathsf{max},0,0]\}$ 
            \hfill ({\em init the heap with the first interval})\label{alg:naive_initheap}
            
        \STATE $\varepsilon_\mathsf{xs} \leftarrow 0$,
            $q_\mathsf{xs} \leftarrow 0$
            \hfill ({\em init the excess error and integral})\label{alg:naive_initxs}

        \WHILE{$\sum_{\varepsilon_i \in H} \varepsilon_i > \tau$} \label{alg:naive_while}
            
            \STATE $k \leftarrow \arg\max_k \varepsilon_k$
                \hfill ({\em get the index of the interval with the largest error})\label{alg:naive_pop}
                
            \STATE $m \leftarrow (a_k+b_k)/2$, $h \leftarrow (b_k-a_k)/2$
            
            \STATE $\mathsf{split} \leftarrow \mathbf{false}$
                \hfill ({\em init $\mathsf{split}$})
            
            \IF{$\varepsilon_k < |q_k|\varepsilon_\mathsf{mach}\mbox{cond}(\mathbf{P}^{(d_k)})\ \vee$ interval too small} \label{alg:naive_drop}
                \STATE $\varepsilon_\mathsf{xs} \leftarrow \varepsilon_\mathsf{xs} + \varepsilon_k$,
                    $q_\mathsf{xs} \leftarrow q_\mathsf{xs} + q_k$
                    \hfill ({\em collect the excess error and integral})\label{alg:naive_addxs}
                \STATE  $H \leftarrow H \setminus \{[a_k,b_k,\mathbf{c}^\mathsf{old},q_k,\varepsilon_k,d_k,\mathsf{nr}_\mathsf{div},\mathsf{nr}_\mathsf{rec}]\}$
                    \hfill ({\em remove the $k$th interval})
            \ELSIF{$d_k < d_\mathsf{max}$} \label{alg:naive_incrdeg}
                \STATE $d_k \leftarrow d_k + 1$
                    \hfill ({\em increase the degree in this interval})\label{alg:naive_update1}\label{alg:naive_fx2}
                \IFOR{$i=0 \dots n_{d_k}$}{$f_i \leftarrow f\left(m - hx^{(d_k)}_i/2\right)$}
                \STATE $\mathbf{c}^{(d_k)} = (\mathbf{P}^{(d_k)})^{-1} \mathbf{f}$
                    \hfill ({\em compute the new interpolation coefficients})
                \STATE $q_k \leftarrow h \boldsymbol\omega^\mathsf{T}\mathbf c^{(d_k)}$
                    \hfill ({\em approximate the new integral})
                \STATE $\varepsilon_k \leftarrow h \|\mathbf{c}^{(d_k)}-\mathbf{c}^\mathsf{old}\|$
                    \hfill ({\em approximate the new error})\label{alg:naive_update2}
                \STATE $\mathsf{split} \leftarrow \frac{\|\mathbf{c}^{(d_k)}-\mathbf{c}^\mathsf{old}\|}{\|\mathbf{c}^{(d_k)}\|} > \mathsf{hint}$
                    \hfill ({\em check change in the coefficients})\label{alg:naive_forcesplit}
            \ELSE
                \STATE $\mathsf{split} \leftarrow \mathbf{true}$
                    \hfill ({\em split the interval if we are already at highest-degree rule\em})
            \ENDIF
            
            \IF{$\mathsf{split} $}
                \STATE  $H \leftarrow H \setminus \{[a_k,b_k,\mathbf{c}^\mathsf{old},q_k,\varepsilon_k,d_k,\mathsf{nr}^\mathsf{old}_\mathsf{div},\mathsf{nr}_\mathsf{rec}]\}$
                    \hfill ({\em remove the $k$th interval})
                \FOR{$i=0 \dots n_{0}$}
                    \STATE $f^\mathsf{left}_i \leftarrow f\left((a_k+m)/2 + hx^{(0)}_i/2\right)$,
                        $f^\mathsf{right}_i \leftarrow f\left((m+b_k)/2 + hx^{(0)}_i/2\right)$\label{alg:naive_fx}
                \ENDFOR
                \FOR{$\mathsf{half} \in \{\mathsf{left},\mathsf{right}\}$}
                    \STATE $\mathbf{c}^\mathsf{half} = (\mathbf{P}^{(0)})^{-1} \mathbf{f}^\mathsf{half}$
                        \hfill ({\em compute the new interpolation coefficients})\label{alg:naive_chalf}
                    \STATE $q_\mathsf{half} \leftarrow h c^\mathsf{half}_0 / \sqrt{2}$
                        \hfill ({\em approximate the new integral})\label{alg:naive_qhalf}
                    \IIF{$q_\mathsf{half} \geq q^{(0)}_k$}{$\mathsf{nr}^\mathsf{half}_\mathsf{div} \leftarrow \mathsf{nr}^\mathsf{old}_\mathsf{div}+1$ {\bf else} $\mathsf{nr}^\mathsf{half}_\mathsf{div} \leftarrow \mathsf{nr}^\mathsf{old}_\mathsf{div}$} \label{alg:naive_ndiv1}
                    \IIF{$\mathsf{nr}^\mathsf{half}_\mathsf{div} > \mathsf{nr}_\mathsf{divmax} \wedge 2\mathsf{nr}^\mathsf{half}_\mathsf{div} > \mathsf{nr}_\mathsf{rec}$}{{\bf return error}}
                        \\ \hfill ({\em abort on divergence})\label{alg:naive_ndiv2}
                    \STATE $\varepsilon_\mathsf{half} \leftarrow h \|\mathbf{c}^\mathsf{half}-\mathbf{T}^{\mathsf{half}}\mathbf{c}^\mathsf{old}\|$
                        \hfill ({\em approximate the new error})\label{alg:naive_ehalf}
                \ENDFOR
                \STATE $H \leftarrow H \cup \{ [a_k,m,\mathbf{c}^\mathsf{left},q_\mathsf{left},\varepsilon_\mathsf{left},0,\mathsf{nr}^\mathsf{left}_\mathsf{div},\mathsf{nr}_\mathsf{rec}+1] , $
                \STATE \hspace{1cm} $ [m,b_k,\mathbf{c}^\mathsf{right},q_\mathsf{right},\varepsilon_\mathsf{right},0,\mathsf{nr}^\mathsf{right}_\mathsf{div},\mathsf{nr}_\mathsf{rec}+1]\}$ 
                    \\ \hfill ({\em push the new intervals back on the heap})
            \ENDIF
            
        \ENDWHILE
        
        \RETURN{$\left[ q_\mathsf{xs} + \sum_{q_i \in H}q_i\right], \left[\varepsilon_\mathsf{xs} + \sum_{\varepsilon_i \in H}\varepsilon_i\right] $}
            \hfill ({\em return the integral and the error})\label{alg:naive_return}
    \end{algorithmic}
\end{algorithm}

The first algorithm (Algorithm~\ref{alg:final_naive}) uses the naive
error estimate in \eqn{err_naive} in a doubly adaptive
strategy using $d_\mathsf{max}+1$ rules of degree $n_i = 2n_{i-1}$,
$i=1 \dots d_\mathsf{max}$, using the nodes $\mathbf{x}^{(i)}$
and the Vandermonde-like matrices $\mathbf{P}^{(i)}$,
$i=0 \dots d_\mathsf{max}$.
For the tests in Section~\ref{sec:validation}, $n_0 = 4$, 
$\mathsf{hint}=0.1$ and $d_\mathsf{max} = 3$ were used.

In Lines~\ref{alg:naive_init1} to \ref{alg:naive_init2}, the coefficients
of the two highest-degree rules are computed and used to approximate
the initial integral and error estimate.
In Line~\ref{alg:naive_initheap} the heap $H$ is initialized with this
interval data.
The algorithm then loops until the sum of the errors over all the intervals
is below the required tolerance (Line~\ref{alg:naive_while}).
At the top of the loop, the interval with the largest error is selected
(Line~\ref{alg:naive_pop}).
If the error in this interval is below the numerical accuracy 
available for the rule used or the interval is too small (\ie the
space between the first two or last two nodes is zero, Line~\ref{alg:naive_drop}),
the interval is dropped and its error and integral are accumulated in the
excess variables $\varepsilon_\mathsf{xs}$ and $q_\mathsf{xs}$
(Line~\ref{alg:naive_addxs}).
If the selected interval has not already used the highest-degree rule
(Line~\ref{alg:naive_incrdeg}), the coefficients of the
next-higher degree rule are computed and the integral and
error estimate are updated (Lines~\ref{alg:naive_update1} to \ref{alg:naive_update2}).
The interval is bisected if either the highest-degree rule
has already been applied or if when increasing the degree
of the rule the coefficients change too much 
(Line~\ref{alg:naive_forcesplit}), analogously to the decision process
suggested by \citeN{ref:Venter2002}.
For the two new sub-intervals, the coefficients for the lowest-degree
rule are computed (Line~\ref{alg:naive_chalf}) and used to approximate
the integral (Line~\ref{alg:naive_qhalf}).
The number of times the integral increases over the sub-interval
is counted in the variables $\mathsf{nr}_\mathsf{div}$ (Line~\ref{alg:naive_ndiv1})
and if they exceed $\mathsf{nr}_\mathsf{divmax}$ and half of the recursion
depth $\mathsf{nr}_\mathsf{rec}$ of that interval, the algorithm aborts 
(Line~\ref{alg:naive_ndiv2}) as per Section~\ref{sec:singularities}.
Note that since the test in \eqn{err_test} requires that both
estimates be of the same degree, we will use, for the estimate $q_k^{(0)}$
from the parent interval, the estimate which was computed for the first rule.
The error estimate for the new interval is computed by transforming
the interpolation coefficients from the parent interval using
\eqn{tlr} and using its difference to the
interpolation in the new interval (Line~\ref{alg:naive_ehalf}).
When the sum of the errors falls below the required tolerance,
the algorithm returns its approximations to the integral and
the integration error (Line~\ref{alg:naive_return}).

\begin{algorithm}
    \caption{int\_refined $(f,a,b,\tau)$}
    \label{alg:final_refined}
    \begin{algorithmic}[1]
    
        \IFOR{$i=0 \dots n$}{$f_i \leftarrow f\left((a+b)/2 - (a-b)x_i/2\right)$}
            \\ \hfill ({\em evaluate the integrand at the nodes $x_i$})\label{alg:refined_init1}
        
        \STATE $\mathbf{c} \leftarrow \mathbf{P}^{-1} \mathbf{f}$,
            \\ \hfill ({\em compute the interpolation coefficients})
            
        \STATE $q_0 \leftarrow (b - a) \boldsymbol\omega^\mathsf{T} \mathbf c^{(d_\mathsf{max})} / 2$
            \hfill ({\em approximate the integral})
        
        \STATE $\varepsilon_0 \leftarrow \infty$
            \hfill ({\em start with a somewhat pessimistic appreciation of the error})\label{alg:refined_initerr}
        
        \STATE $H \leftarrow \{[a,b,\mathbf{c},\mathbf{b},q_0,\varepsilon_0,0,0]\}$ 
            \hfill ({\em init the heap with the first interval})\label{alg:refined_init2}
            
        \STATE $\varepsilon_\mathsf{xs} \leftarrow 0$,
            $q_\mathsf{xs} \leftarrow 0$
            \hfill ({\em init the excess error and integral})

        \WHILE{$\sum_{\varepsilon_i \in H} \varepsilon_i > \tau$} \label{alg:refined_loop}
            
            \STATE $k \leftarrow \arg\max_k \varepsilon_k$
                \hfill ({\em get the index of the interval with the largest error})\label{alg:refined_pop}
                
            \STATE  $H \leftarrow H \setminus \{[a_k,b_k,\mathbf{c}^\mathsf{old},\mathbf{b}^\mathsf{old},q_k,\varepsilon_k,\mathsf{nr}^\mathsf{old}_\mathsf{div},\mathsf{nr}_\mathsf{rec}]\}$
                \hfill ({\em remove the $k$th interval})\label{alg:refined_pop2}

            \IF{$\varepsilon_k < |q_k|\varepsilon_\mathsf{mach}\mbox{cond}(\mathbf{P})\ \vee$ interval too small} \label{alg:refined_drop}
                \STATE $\varepsilon_\mathsf{xs} \leftarrow \varepsilon_\mathsf{xs} + \varepsilon_k$,
                    $q_\mathsf{xs} \leftarrow q_\mathsf{xs} + q_k$
                    \hfill ({\em collect the excess error and integral})\label{alg:refined_xs}
            \ELSE
                \STATE $m \leftarrow (a_k+b_k)/2$, $h \leftarrow (b_k-a_k)/2$
                \FOR{$i=0 \dots n$}
                    \STATE $f^\mathsf{left}_i \leftarrow f\left((a_k+m)/2 + hx_i/2\right)$,
                        $f^\mathsf{right}_i \leftarrow f\left((m+b_k)/2 + hx_i/2\right)$
                        \\ \hfill ({\em evaluate the integrand at the nodes $x_i$ in the sub-intervals})\label{alg:refined_fx}
                \ENDFOR
                \FOR{$\mathsf{half} \in \{\mathsf{left},\mathsf{right}\}$}
                    \STATE $\mathbf{c}^\mathsf{half} = (\mathbf{P})^{-1} \mathbf{f}^\mathsf{half}$
                        \hfill ({\em compute the new interpolation coefficients})\label{alg:refined_chalf}
                    \STATE $q_\mathsf{half} \leftarrow h \boldsymbol\omega^\mathsf{T} \mathbf c^\mathsf{half}$
                        \hfill ({\em approximate the new integral})\label{alg:refined_qhalf}
                    \IIF{$q_\mathsf{half} \geq q_k$}{$\mathsf{nr}^\mathsf{half}_\mathsf{div} \leftarrow \mathsf{nr}^\mathsf{old}_\mathsf{div}+1$ {\bf else} $\mathsf{nr}^\mathsf{half}_\mathsf{div} \leftarrow \mathsf{nr}^\mathsf{old}_\mathsf{div}$} \label{alg:refined_ndiv}
                    \IIF{$\mathsf{nr}^\mathsf{half}_\mathsf{div} > \mathsf{nr}_\mathsf{divmax} \wedge 2\mathsf{nr}^\mathsf{half}_\mathsf{div} > \mathsf{nr}_\mathsf{rec}$}{{\bf return error}}
                        \\ \hfill ({\em abort on divergence})\label{alg:refined_div}
                    \STATE $f^{(n+1)}_\mathsf{half} \leftarrow \frac{\|\mathbf{c}^\mathsf{half}-\mathbf{T}^{(\mathsf{half})}\mathbf{c}^\mathsf{old}\|}{\|\mathbf{b}^\mathsf{half}-2^{n+1}\mathbf{T}^{(\mathsf{half})}\mathbf{b}^\mathsf{old}\|}$
                        \hfill ({\em approximate the higher derivative})
                    \IF{$\max \left\{ \left|\mathbf{P}\mathbf{c}^\mathsf{old} - \mathbf{f}^\mathsf{half}\right| - \vartheta_1 f^{(n+1)}_\mathsf{half} \left| \mathbf{P}\mathbf{b}^\mathsf{old} \right|\right\} > 0$} \label{alg:refined_test}
                        \STATE $\varepsilon_\mathsf{half} \leftarrow h \|\mathbf{c}^\mathsf{half} - \mathbf{c}^\mathsf{old}\|$
                            \hfill ({\em compute the un-scaled error})\label{alg:refined_unscaled}
                    \ELSE
                        \STATE $\varepsilon_\mathsf{half} \leftarrow h f^{(n+1)}_\mathsf{half} \|\mathbf{b}^\mathsf{half}\|$
                            \hfill ({\em compute the extrapolated error})\label{alg:refined_scaled}
                    \ENDIF
                \ENDFOR
                \STATE $H \leftarrow H \cup \{ [a_k,m,\mathbf{c}^\mathsf{left},\mathbf{b}^\mathsf{left},q_\mathsf{left},\varepsilon_\mathsf{left},\mathsf{nr}^\mathsf{left}_\mathsf{div},\mathsf{nr}_\mathsf{rec}+1] $, \dots \\
                    \hspace{1cm} $[m,b_k,\mathbf{c}^\mathsf{right},\mathbf{b}^\mathsf{right},q_\mathsf{right},\varepsilon_\mathsf{right},\mathsf{nr}^\mathsf{left}_\mathsf{div},\mathsf{nr}_\mathsf{rec}+1]\}$ 
                    \\ \hfill ({\em push the new intervals back on the heap})\label{alg:refined_push}
            \ENDIF
            
        \ENDWHILE
        
        \RETURN{$\left[ q_\mathsf{xs} + \sum_{q_i \in H}q_i\right], \left[ \varepsilon_\mathsf{xs} + \sum_{\varepsilon_i \in H}\varepsilon_i\right]$}
            \hfill ({\em return the integral and the error})\label{alg:refined_return}
            
    \end{algorithmic}
\end{algorithm}

The second algorithm (Algorithm~\ref{alg:final_refined}) uses the refined
error estimate in \eqn{err_refined}, which re-uses
the coefficients from a previous level of recursion.
For the results in Section~\ref{sec:validation}, $n=10$ and
$\vartheta_1 = 1.1$ were used.

In Lines~\ref{alg:refined_init1} to \ref{alg:refined_init2} an initial estimate
is computed and used to initialize the heap $H$.
The error estimate is set to $\infty$ since it can not be estimated
(Line~\ref{alg:refined_initerr}).
While the sum of error estimates is above the required tolerance,
the algorithm selects the interval with the largest error estimate
(Line~\ref{alg:refined_pop}) and removes it from the heap
(Line~\ref{alg:refined_pop2}).
As with the previous algorithm, if the error estimate is smaller
than the numerical precision of the integral or the interval is too
small, the interval is dropped (Line~\ref{alg:refined_drop})
and its error and integral estimates are stored in the excess variables
$\varepsilon_\mathsf{xs}$ and $q_\mathsf{xs}$ (Line~\ref{alg:refined_xs}).
The algorithm then computes the new coefficients for each
sub-interval (Line~\ref{alg:refined_chalf}), as well
as their integral approximation (Line~\ref{alg:refined_qhalf}).
If the integral over the sub-interval is larger than over the
previous interval, the variable $\mathsf{nr}_\mathsf{div}$ is increased
(Line~\ref{alg:refined_ndiv}) and if it exceeds $\mathsf{nr}_\mathsf{divmax}$ and
half of the recursion depth $\mathsf{nr}_\mathsf{rec}$,
the algorithm aborts with an error (Line~\ref{alg:refined_div}).
In Line~\ref{alg:refined_test} the algorithm tests whether the
conditions laid out in \eqn{err_test}
for the approximation of the $n+1$st derivative hold.
If they do not, the un-scaled error estimate is returned
(Line~\ref{alg:refined_unscaled}), otherwise, the scaled
estimate is returned (Line~\ref{alg:refined_scaled}).
Finally, both sub-intervals are returned to the heap
(Line~\ref{alg:refined_push}).
Once the required tolerance is met, the algorithm returns its 
approximations to the integral and
the integration error (Line~\ref{alg:refined_return}).

Although the algorithm descriptions in Algorithms~\ref{alg:final_naive}
and \ref{alg:final_refined} are quite complete, some details
have been omitted for simplicity.
First of all, when the function values are computed
(Lines~\ref{alg:naive_init1}, \ref{alg:naive_fx2} and \ref{alg:naive_fx} of 
Algorithm~\ref{alg:final_naive} and Lines~\ref{alg:refined_init1} and
\ref{alg:refined_fx} of Algorithm~\ref{alg:final_refined}), it is
understood that previously computed function values at the same
nodes, \ie on the edges of the domain for both algorithms or
inside the Clenshaw-Curtis rules of increasing degree for
Algorithm~\ref{alg:final_naive}, are re-used and not re-evaluated.

We have also not included the downdate of the interpolations when
$\mathsf{NaN}$ or $\pm\mathsf{Inf}$ is encountered.
This is done as is shown in Algorithm~\ref{alg:sing}.
If, when evaluating the integrand, a non-numerical value
is encountered, the function values is set to zero and
the index of the node is stored in $\mathsf{nans}$
(Line~\ref{alg:sing_catch}).
The coefficients of the interpolation are then computed
for those function values (Line~\ref{alg:sing_c}).
For each index in $\mathsf{nans}$, first the coefficients $\mathbf{b}$
of the Newton polynomial over the nodes of the quadrature rule
are down-dated as per \eqn{sing_downdateb} (Line~\ref{alg:sing_b}).
The downdated $\mathbf{b}$ is then in turn used to downdate
the interpolation coefficients $\mathbf{c}$ as per
\eqn{sing_downdatec} (Line~\ref{alg:sing_c2}).

Furthermore, to improve memory efficiency, both algorithms
maintain at most 200 intervals in the heap.
If this number is exceeded, the interval with the smallest
error estimate is removed and its integral and error estimates
are added to the excess variables $q_\mathsf{xs}$ and
$\varepsilon_\mathsf{xs}$ respectively.

\begin{algorithm}
    \caption{Interpolation downdate procedure}
    \label{alg:sing}
    \begin{algorithmic}[1]
        \STATE $\mathsf{nans} \leftarrow \{\}$
            \hfill ({\em initialize $\mathsf{nans}$})
        \FOR{$i=0 \dots n$}
            \STATE $f_i \leftarrow f\left(a + \frac{x_i + 1}{2}(b-a)\right)$
                \hfill ({\em evaluate the integrand at the nodes $x_i$})
            \IIF{$f_i \in \{\mathsf{NaN},\mathsf{Inf}\}$}{$f_i \leftarrow 0, \mathsf{nans} \leftarrow \mathsf{nans} \cup \{i\}$}
                \\ \hfill ({\em if the result is non-numerical, set the node to zero and remember it})\label{alg:sing_catch}
        \ENDFOR
        \STATE $\mathbf{c} \leftarrow \mathbf{P}^{-1} \mathbf{f}$ \label{alg:sing_c}
            \hfill ({\em compute the initial interpolation coefficients})
        \FOR{$i \in \mathsf{nans}$}
            \STATE $\mathbf{b} \leftarrow \mathbf{U}^{-1}_i \mathbf{b}$
                \hfill ({\em downdate the coefficients of the Newton polynomial})\label{alg:sing_b}
            \STATE $\mathbf{c} \leftarrow \mathbf{c} - \frac{c_n}{b_n}\mathbf{b}$
                \hfill ({\em downdate the coefficients of the interpolation})\label{alg:sing_c2}
            \STATE $n \leftarrow n - 1$ \label{alg:final1_downdate}
                \hfill ({\em decrement the degree})
        \ENDFOR
    \end{algorithmic}
\end{algorithm}

\section{Validation}
\label{sec:validation}

\begin{sidewaystable}
    \begin{center}\begin{scriptsize}
    \begin{tabular}{l|c|c|c|c|c|c|c|c|c|c|c|c|c|c|c}
        $\tau=10^{-3}$ & \multicolumn{3}{c|}{\tt quadl} & \multicolumn{3}{c|}{\tt DQAGS} & \multicolumn{3}{c|}{\tt da2glob} & \multicolumn{3}{c|}{Algorithm~\ref{alg:final_refined}} & \multicolumn{3}{c}{Algorithm~\ref{alg:final_naive}} \\
        $f(x)$ & \Checkmark & \XSolid & $n_\mathsf{eval}$ & \Checkmark & \XSolid & $n_\mathsf{eval}$ & \Checkmark & \XSolid & $n_\mathsf{eval}$ & \Checkmark & \XSolid & $n_\mathsf{eval}$ & \Checkmark & \XSolid & $n_\mathsf{eval}$ \\ \hline
Eqn (\ref{eqn:val_sing})  & $395$ & $605$ & $94.47$ & $915$ & $85$ & $450.62$ & $977$ & $23$ & $95.97$ & $1000$ & $0$ & $361.63$ & $1000$ & $0$ & $281.75$ \\ \hline
Eqn (\ref{eqn:val_disc})  & $885$ & $115$ & $117.48$ & $963$ & $37$ & $399.13$ & $1000$ & $0$ & $59.62$ & $1000$ & $0$ & $306.80$ & $1000$ & $0$ & $175.30$ \\ \hline
Eqn (\ref{eqn:val_c0})    & $867$ & $133$ & $46.14$ & $1000$ & $0$ & $179.13$ & $1000$ & $0$ & $29.14$ & $1000$ & $0$ & $99.39$ & $1000$ & $0$ & $113.36$ \\ \hline
Eqn (\ref{eqn:val_peak})  & $360$ & $640$ & $88.53$ & $706$ & $294\,(16)$ & $409.37$ & $882$ & $118.5$ & $136.50$ & $1000$ & $0$ & $498.56$ & $1000$ & $0$ & $342.02$ \\ \hline
Eqn (\ref{eqn:val_4peak}) & $265$ & $735$ & $291.60$ & $962$ & $38\,(4)$ & $1715.74$ & $997$ & $3$ & $480.64$ & $1000$ & $0$ & $1457.63$ & $996$ & $4$ & $990.57$ \\ \hline
Eqn (\ref{eqn:val_oscill})& $995$ & $5$ & $765.39$ & $1000$ & $0$ & $450.53$ & $1000$ & $0$ & $523.35$ & $1000$ & $0$ & $688.98$ & $1000$ & $0$ & $879.37$ \\ \hline
    \end{tabular}\end{scriptsize}\end{center}
%
    \begin{center}\begin{scriptsize}\begin{tabular}{l|c|c|c|c|c|c|c|c|c|c|c|c|c|c|c}
        $\tau=10^{-6}$ & \multicolumn{3}{c|}{\tt quadl} & \multicolumn{3}{c|}{\tt DQAGS} & \multicolumn{3}{c|}{\tt da2glob} & \multicolumn{3}{c|}{Algorithm~\ref{alg:final_refined}} & \multicolumn{3}{c}{Algorithm~\ref{alg:final_naive}} \\
        $f(x)$ & \Checkmark & \XSolid & $n_\mathsf{eval}$ & \Checkmark & \XSolid & $n_\mathsf{eval}$ & \Checkmark & \XSolid & $n_\mathsf{eval}$ & \Checkmark & \XSolid & $n_\mathsf{eval}$ & \Checkmark & \XSolid & $n_\mathsf{eval}$ \\ \hline
Eqn (\ref{eqn:val_sing})  & $406$ & $594$ & $371.64$ & $910$ & $90$ & $1094.27$ & $993\,(7)$ & $7$ & $298.45$ & $1000$ & $0$ & $993.84$ & $1000$ & $0$ & $870.30$ \\ \hline
Eqn (\ref{eqn:val_disc})  & $879$ & $121$ & $237.54$ & $906$ & $94$ & $791.28$ & $1000$ & $0$ & $100.12$ & $1000$ & $0$ & $626.96$ & $1000$ & $0$ & $316.15$ \\ \hline
Eqn (\ref{eqn:val_c0})    & $780$ & $220$ & $103.56$ & $987$ & $13$ & $365.32$ & $1000$ & $0$ & $58.12$ & $1000$ & $0$ & $255.25$ & $1000$ & $0$ & $313.78$ \\ \hline
Eqn (\ref{eqn:val_peak})  & $987$ & $13$ & $482.85$ & $931$ & $69\,(69)$ & $699.22$ & $1000$ & $0$ & $290.52$ & $1000$ & $0$ & $766.88$ & $1000$ & $0$ & $616.96$ \\ \hline
Eqn (\ref{eqn:val_4peak}) & $994$ & $6$ & $1318.02$ & $999$ & $1$ & $2075.72$ & $1000$ & $0$ & $870.0$ & $1000$ & $0$ & $2292.21$ & $1000$ & $0$ & $1840.57$ \\ \hline
Eqn (\ref{eqn:val_oscill})& $999$ & $1$ & $2077.80$ & $1000$ & $0$ & $581.41$ & $1000\,(1)$ & $0$ & $695.02$ & $1000$ & $0$ & $1193.97$ & $1000$ & $0$ & $1200.88$ \\ \hline
    \end{tabular}\end{scriptsize}\end{center}
%
    \begin{center}\begin{scriptsize}\begin{tabular}{l|c|c|c|c|c|c|c|c|c|c|c|c|c|c|c}
        $\tau=10^{-9}$ & \multicolumn{3}{c|}{\tt quadl} & \multicolumn{3}{c|}{\tt DQAGS} & \multicolumn{3}{c|}{\tt da2glob} & \multicolumn{3}{c|}{Algorithm~\ref{alg:final_refined}} & \multicolumn{3}{c}{Algorithm~\ref{alg:final_naive}} \\
        $f(x)$ & \Checkmark & \XSolid & $n_\mathsf{eval}$ & \Checkmark & \XSolid & $n_\mathsf{eval}$ & \Checkmark & \XSolid & $n_\mathsf{eval}$ & \Checkmark & \XSolid & $n_\mathsf{eval}$ & \Checkmark & \XSolid & $n_\mathsf{eval}$ \\ \hline
Eqn (\ref{eqn:val_sing})  & $360$ & $640\,(38)$ & $1064.38$ & $394\,(389)$ & $606\,(597)$ & $1525.27$ & $874\,(222)$ & $126\,(125)$ & $627.40$ & $897\,(163)$ & $103\,(103)$ & $2007.81$ & $894\,(23)$ & $106\,(99)$ & $1852.39$ \\ \hline
Eqn (\ref{eqn:val_disc})  & $877$ & $123$ & $357.51$ & $609\,(571)$ & $391\,(284)$ & $1025.39$ & $1000$ & $0$ & $141.82$ & $1000$ & $0$ & $945.20$ & $1000$ & $0$ & $460.87$ \\ \hline
Eqn (\ref{eqn:val_c0})    & $779$ & $221$ & $186.24$ & $867\,(734)$ & $133\,(98)$ & $485.73$ & $1000$ & $0$ & $93.74$ & $1000$ & $0$ & $415.14$ & $1000$ & $0$ & $522.28$ \\ \hline
Eqn (\ref{eqn:val_peak})  & $998$ & $2$ & $1218.57$ & $931\,(414)$ & $69\,(69)$ & $805.81$ & $1000$ & $0$ & $474.81$ & $1000$ & $0$ & $1291.34$ & $1000$ & $0$ & $1078.27$ \\ \hline
Eqn (\ref{eqn:val_4peak}) & $1000$ & $0$ & $3375.24$ & $999\,(936)$ & $1\,(1)$ & $2369.30$ & $1000$ & $0$ & $1431.70$ & $1000$ & $0$ & $3915.23$ & $1000$ & $0$ & $3292.26$ \\ \hline
Eqn (\ref{eqn:val_oscill})& $945$ & $55\,(55)$ & $5274.96$ & $1000\,(873)$ & $0$ & $672.17$ & $998\,(63)$ & $2\,(2)$ & $799.86$ & $1000\,(3)$ & $0$ & $2028.08$ & $1000\,(1)$ & $0$ & $1344.78$ \\ \hline
    \end{tabular}\end{scriptsize}\end{center}
%
    \begin{center}\begin{scriptsize}\begin{tabular}{l|c|c|c|c|c|c|c|c|c|c|c|c|c|c|c}
        $\tau=10^{-12}$ & \multicolumn{3}{c|}{\tt quadl} & \multicolumn{3}{c|}{\tt DQAGS} & \multicolumn{3}{c|}{\tt da2glob} & \multicolumn{3}{c|}{Algorithm~\ref{alg:final_refined}} & \multicolumn{3}{c}{Algorithm~\ref{alg:final_naive}} \\
        $f(x)$ & \Checkmark & \XSolid & $n_\mathsf{eval}$ & \Checkmark & \XSolid & $n_\mathsf{eval}$ & \Checkmark & \XSolid & $n_\mathsf{eval}$ & \Checkmark & \XSolid & $n_\mathsf{eval}$ & \Checkmark & \XSolid & $n_\mathsf{eval}$ \\ \hline
Eqn (\ref{eqn:val_sing})  & $265$ & $735\,(359)$ & $2868.37$ & $0$ & $1000\,(1000)$ & $1525.27$ & $494\,(267)$ & $506\,(506)$ & $1068.39$ & $558\,(187)$ & $442\,(442)$ & $9156.03$ & $551\,(65)$ & $449\,(449)$ & $8713.21$ \\ \hline
Eqn (\ref{eqn:val_disc})  & $892$ & $108$ & $489.84$ & $3\,(3)$ & $997\,(900)$ & $1025.39$ & $1000$ & $0$ & $184.42$ & $1000\,(2)$ & $0$ & $1264.45$ & $1000$ & $0$ & $606.67$ \\ \hline
Eqn (\ref{eqn:val_c0})    & $781$ & $219$ & $317.73$ & $19\,(15)$ & $981\,(949)$ & $485.73$ & $1000$ & $0$ & $133.60$ & $1000$ & $0$ & $575.06$ & $1000$ & $0$ & $737.59$ \\ \hline
Eqn (\ref{eqn:val_peak})  & $942$ & $58$ & $3235.86$ & $871\,(851)$ & $129\,(127)$ & $845.92$ & $996\,(765)$ & $4\,(4)$ & $710.62$ & $1000\,(592)$ & $0$ & $10626.40$ & $1000\,(428)$ & $0$ & $17855.83$ \\ \hline
Eqn (\ref{eqn:val_4peak}) & $711$ & $289\,(289)$ & $8799.54$ & $987\,(987)$ & $13\,(13)$ & $2400.09$ & $1000\,(849)$ & $0$ & $2068.38$ & $1000\,(410)$ & $0$ & $9619.10$ & $1000\,(352)$ & $0$ & $13596.62$ \\ \hline
Eqn (\ref{eqn:val_oscill})& $223$ & $777\,(777)$ & $9640.14$ & $951\,(951)$ & $49\,(49)$ & $672.17$ & $548\,(548)$ & $452\,(452)$ & $894.86$ & $989\,(658)$ & $11\,(11)$ & $9808.66$ & $987\,(514)$ & $13\,(13)$ & $22433.79$ \\ \hline
    \end{tabular}\end{scriptsize}\end{center}
    \caption{Results of the Lyness-Kaganove tests for $\tau=10^{-3}, 10^{-6}, 10^{-9}$
        and $10^{-12}$. The columns marked with \Checkmark and \XSolid indicate the number
        of correct and incorrect results respectively, out of 1\,000 runs. The numbers
        in brackets indicate the number of runs in which a warning was issued. The column
        $n_\mathsf{eval}$ contains the average number of function evaluations required
        for each run.}
    \label{tab:LKtest}
\end{sidewaystable}

In the following, we will test both algorithms described in
Section~\ref{sec:algorithm} against the following routines:
\begin{itemize}
    
    \item {\tt quadl}, {\small MATLAB}'s adaptive quadrature routine
        \cite{ref:Mathworks2003}, based on Gander and Gautschi's 
        {\tt adaptlob} \cite{ref:Gander2001}.
        This algorithm uses a 4-point Gauss-Lobatto quadrature rule
        and its 7-point Kronrod extension.
                
    \item {\tt DQAGS}, {\small QUADPACK}'s \cite{ref:Piessens1983} adaptive quadrature
        routine using a 10-point Gauss quadrature rule and its
        21-point Kronrod extension as well as the $\varepsilon$-Algorithm
        to extrapolate the integral and error estimate.
        The routine is called through the GNU Octave \cite{ref:Eaton2002}
        package's {\tt quad} routine.
        
    \item {\tt da2glob} from \citeN{ref:Espelid2007}, which uses a
        doubly-adaptive strategy over rules of degree
        5, 9, 17 and 27 over 5, 9, 17 and 33 equidistant nodes\footnote{
        The final rule of degree 27 over 33 nodes is constructed such
        as to maximize its numerical stability.} respectively,
        using the error estimator described in
        \citeN{ref:Berntsen1991}.

\end{itemize}

The two algorithms in Section~\ref{sec:algorithm} were implemented in
the {\small MATLAB} programming language\footnote{The source-code of both
routines is available online at {\tt http://people.inf.ethz.ch/gonnetp/toms/}}.

Over the years, several authors have specified sets of test functions
to evaluate the performance and reliability of quadrature routines.
In the following, we will use, with some minor modifications,
the test ``families'' suggested by
\citeN{ref:Lyness1977} and the ``battery'' of functions
compiled by \citeN{ref:Gander1998}, which are
an extension of the set proposed by \citeN{ref:Kahaner1971}.

The function families used for the Lyness-Kaganove test are
\begin{eqnarray}
    \int_0^1 |x-\lambda|^\alpha\,\mbox{d}x,\quad & & \lambda \in [0,1],\ \alpha \in [-0.5,0] \label{eqn:val_sing} \\
    \int_0^1 (x>\lambda) e^{\alpha x}\,\mbox{d}x,\quad & & \lambda \in [0,1],\ \alpha \in [0,1] \label{eqn:val_disc} \\
    \int_0^1 \exp (-\alpha|x-\lambda|)\,\mbox{d}x,\quad & & \lambda \in [0,1],\ \alpha \in [0,4] \label{eqn:val_c0} \\
    \int_1^2 10^\alpha / ((x-\lambda)^2 + 10^\alpha)\,\mbox{d}x,\quad & & \lambda \in [1,2], \alpha \in [-6,-3] \label{eqn:val_peak} \\
    \int_1^2 \sum_{i=1}^4 10^\alpha / ((x-\lambda_i)^2 + 10^\alpha)\,\mbox{d}x,\quad & & \lambda_i \in [1,2], \alpha \in [-5,-3] \label{eqn:val_4peak} \\
    \int_0^1 2\beta(x-\lambda)\cos(\beta(x-\lambda)^2)\,\mbox{d}x, \quad & & \lambda \in [0,1],\ \alpha \in [1.8,2], \label{eqn:val_oscill} \\
    & & \beta = 10^\alpha/\max\{\lambda^2,(1-\lambda)^2\} \nonumber
\end{eqnarray}
where the boolean expressions are evaluated to 0 or 1.
The integrals are computed to relative precisions\footnote{
Since all algorithms tested use an absolute error bound, the
exact integral times the relative tolerance was used.} of $\tau=10^{-3}$,
$10^{-6}$, $10^{-9}$ and $10^{-12}$ for $1\,000$ realizations
of the random parameters $\lambda$ and $\alpha$.
The results of these tests are shown in Table~\ref{tab:LKtest}.
For each function, the number of correct and incorrect integrations
is given with, in brackets, the number of cases each where a warning
(either explicit or whenever an error estimate larger than the requested
tolerance is returned) was issued.
We consider an integration to be correct only when the returned value
is within the required tolerance of the exact result.

The functions used for the ``battery'' test are
\begin{equation*}\begin{array}{rclrcl}
    f_1 & = & \textstyle \int_0^1 e^{x} \,\mbox{d}x &
    
        f_{14} & = & \textstyle \int_0^{10} \sqrt{50}e^{-50\pi x^2} \,\mbox{d}x \\
        
    f_2 & = & \textstyle \int_0^1 (x>0.3) \,\mbox{d}x &
    
        f_{15} & = & \textstyle \int_0^{10} 25 e^{-25x} \,\mbox{d}x \\
        
    f_3 & = & \textstyle \int_0^1 x^{1/2} \,\mbox{d}x &
    
        f_{16} & = & \textstyle \int_0^{10} 50(\pi(2500x^2+1))^{-1} \,\mbox{d}x \\
        
    f_4 & = & \textstyle \int_{-1}^1 (\frac{23}{25}\cosh(x) - \cos(x)) \,\mbox{d}x \quad &
    
        f_{17} & = & \textstyle \int_{0}^1 50(\sin(50\pi x)/(50\pi x))^2 \,\mbox{d}x \\
        
    f_5 & = & \textstyle \int_{-1}^1 (x^4 + x^2 + 0.9)^{-1} \,\mbox{d}x &
    
        f_{18} & = & \textstyle \int_{0}^\pi \cos(\cos(x) + 3 \sin(x) + 2 \cos(2 x) + 3 \cos(3 x)) \,\mbox{d}x \\
        
    f_6 & = & \textstyle \int_0^1 x^{3/2} \,\mbox{d}x &
    
        f_{19} & = & \textstyle \int_0^1 \log(x) \,\mbox{d}x \\
        
    f_7 & = & \textstyle \int_0^1 x^{-1/2} \,\mbox{d}x &
    
        f_{20} & = & \textstyle \int_{-1}^1 (1.005 + x^2)^{-1} \,\mbox{d}x \\
        
    f_8 & = & \textstyle \int_0^1 (1+x^4)^{-1} \,\mbox{d}x &
    
        f_{21} & = & \textstyle \int_0^1 \sum_{i=1}^3 \left[ \cosh(20^i(x-2i/10))\right]^{-1} \,\mbox{d}x \\
        
    f_9 & = & \textstyle \int_0^1 2(2 + \sin(10\pi x))^{-1} \,\mbox{d}x &
    
        f_{22} & = & \textstyle \int_0^1 4\pi^2x \sin(20\pi x) \cos(2\pi x) \,\mbox{d}x \\
        
    f_{10} & = & \textstyle \int_0^1 (1+x)^{-1} \,\mbox{d}x &
    
        f_{23} & = & \textstyle \int_0^1 (1+(230x-30)^2)^{-1} \,\mbox{d}x \\
        
    f_{11} & = & \textstyle \int_0^1 (1+e^x)^{-1} \,\mbox{d}x & 
    
        f_{24} & = & \textstyle \int_0^3 \lfloor e^x \rfloor \,\mbox{d}x \\
        
    f_{12} & = & \textstyle \int_0^1 x(e^x-1)^{-1} \,\mbox{d}x &
    
        f_{25} & = & \textstyle \int_0^5 (x+1)(x<1) + (3-x)(1 \leq x \leq3) \\
        
    f_{13} & = & \textstyle \int_{0}^1 \sin(100\pi x) / (\pi x) \,\mbox{d}x &
        & &  \quad + 2(x>3) \,\mbox{d}x \\
        
\end{array}\end{equation*}
where the boolean expressions
in $f_2$ and $f_{25}$ evaluate to 0 or 1.
The functions are taken from \citeN{ref:Gander1998} with the following
modifications:
\begin{itemize}
    \item No special treatment is given to the case $x=0$ in
        $f_{12}$, allowing the integrand to return $\mathsf{NaN}$.
    \item $f_{13}$ and $f_{17}$ are integrated from 0 to 1 as
        opposed to 0.1 to 1 and 0.01 to 1 respectively, allowing
        the integrand to return $\mathsf{NaN}$ for $x=0$.
    \item No special treatment of $x<10^{-15}$ in $f_{19}$ allowing
        the integrand to return $-\mathsf{Inf}$.
    \item $f_{24}$ was suggested by J.\ Waldvogel as a simple yet tricky test
        function with multiple discontinuities.
    \item $f_{25}$ was introduced in \citeN{ref:Gander1998}, yet
        not used in the battery test therein.
\end{itemize}
The rationale for the modifications of $f_{12}$, $f_{13}$, $f_{17}$ and
$f_{19}$ is that we can't, on one hand, assume that the user was
careful enough or knew enough about his or her integrand
to remove the non-numerical values,
and on the other hand assume that he or she would still resort to
a general-purpose quadrature routine to integrate it.
Any general purpose quadrature routine should be robust
enough to deal with any function, provided by either careful or 
careless users.

These changes have little effect on {\tt quadl} and {\tt DQAGS}
since, as mentioned in Section~\ref{sec:singularities},
the former shifts the integration boundaries by $\varepsilon_\mathsf{mach}$
if a non-numerical value is encountered on the edges of the integration
domain and the later uses Gauss and Gauss-Kronrod quadrature rules
which do not contain the end-points and thus avoid the $\mathsf{NaN}$
returned at $x=0$ for $f_{12}$, $f_{13}$ and $f_{17}$ and the
$\mathsf{Inf}$ at $x=0$ in $f_{19}$.
{\tt da2glob} treats $\mathsf{NaN}$s by setting the offending function
values to 1.
Due to the rather fortunate coincidence that $f_{12}$, $f_{13}$ and
$f_{17}$ all have a limit of 1 for $x\rightarrow 0$, this integrator
is in no way troubled by these integrands.

\begin{sidewaystable}
    \begin{center}\begin{scriptsize}
    \begin{tabular}{l|c|c|c|c|c||c|c|c|c|c||c|c|c|c|c||c|c|c|c|c}
        & \multicolumn{5}{c||}{$\tau = 10^{-3}$}& \multicolumn{5}{c||}{$\tau = 10^{-6}$}& \multicolumn{5}{c||}{$\tau = 10^{-9}$}& \multicolumn{5}{c}{$\tau = 10^{-12}$} \\
        $f(x)$ & {\tt quadl} & {\tt DQAGS} & {\tt da2glob} & Alg.~\ref{alg:final_refined} & Alg.~\ref{alg:final_naive} 
        & {\tt quadl} & {\tt DQAGS} & {\tt da2glob} & Alg.~\ref{alg:final_refined} & Alg.~\ref{alg:final_naive} 
        & {\tt quadl} & {\tt DQAGS} & {\tt da2glob} & Alg.~\ref{alg:final_refined} & Alg.~\ref{alg:final_naive} 
        & {\tt quadl} & {\tt DQAGS} & {\tt da2glob} & Alg.~\ref{alg:final_refined} & Alg.~\ref{alg:final_naive} \\ \hline
$f_{1 }$  & 18 & 21 & {\bf 9} & 27 & 33 & 18 & 21 & {\bf 9} & 27 & 33 & 18 & 21 & {\bf 9} & 27 & 33 & 18 & 21 & {\bf 17} & 27 & 33 \\
$f_{2 }$  & 108 & 357 & {\bf 61} & 283 & 161 & 198 & 357 & {\bf 101} & 603 & 301 & 318 & 357 & {\bf 141} & 923 & 441 & 408 & 357 & {\bf 181} & 1243 & 581 \\
$f_{3 }$  & 48 & 105 & {\bf 25} & 107 & 101 & 108 & 231 & {\bf 65} & 315 & 429 & 258 & 231 & {\bf 137} & 523 & 799 & \st{618} & {\bf 231} & 241 & 827 & 1191 \\
$f_{4 }$  & 18 & 21 & {\bf 9} & 27 & 33 & 18 & 21 & {\bf 9} & 27 & 33 & {\bf 18} & 21 & 25 & 27 & 33 & 48 & {\bf 21} & 33 & 27 & 33 \\
$f_{5 }$  & 18 & 21 & {\bf 17} & 27 & 33 & 48 & {\bf 21} & 33 & 59 & 95 & {\bf 48} & 63 & 65 & 59 & 95 & 168 & {\bf 63} & 65 & 123 & 219 \\
$f_{6 }$  & 18 & 21 & {\bf 9} & 27 & 33 & 48 & 105 & {\bf 41} & 123 & 159 & 108 & 189 & {\bf 73} & 251 & 359 & 288 & 189 & {\bf 137} & 427 & 607 \\
$f_{7 }$  & 289 & 231 & {\bf 121} & 411 & 269 & 439 & {\bf 231} & 285 & 1035 & 709 & \st{889} & {\bf 231} & 581 & 1867 & 1409 & \st{2429} & {\bf 231} & 965 & 3291 & 2179 \\
$f_{8 }$  & 18 & 21 & {\bf 17} & 27 & 33 & {\bf 18} & 21 & 25 & 43 & 33 & 48 & {\bf 21} & 33 & 59 & 95 & 138 & 63 & {\bf 49} & 123 & 95 \\
$f_{9 }$  & 198 & 315 & {\bf 121} & 251 & 261 & 468 & 399 & {\bf 233} & 411 & 587 & 1038 & 567 & {\bf 401} & 891 & 991 & 2808 & 735 & {\bf 577} & 1387 & 1425 \\
$f_{10}$  & 18 & 21 & {\bf 9} & 27 & 33 & 18 & 21 & {\bf 17} & 27 & 33 & 48 & 21 & {\bf 17} & 59 & 33 & 48 & {\bf 21} & 33 & 75 & 33 \\
$f_{11}$  & 18 & 21 & {\bf 9} & 27 & 33 & 18 & 21 & {\bf 9} & 27 & 33 & 18 & 21 & {\bf 9} & 27 & 33 & 48 & 21 & {\bf 17} & 43 & 33 \\
$f_{12}$  & 19 & 21 & {\bf 9} & 59 & 47 & 19 & 21 & {\bf 9} & 59 & 55 & 19 & 21 & {\bf 9} & 59 & 63 & 19 & 21 & {\bf 9} & 59 & 63 \\
$f_{13}$  & \st{589} & {\bf 651} & 929 & 907 & 1403 & 1519 & {\bf 1323} & 1469 & 1643 & 2347 & 4879 & {\bf 1323} & 1913 & 2955 & 2459 & \st{10039} & {\bf 1323} & 2233 & 5035 & 2521 \\
$f_{14}$  & 78 & 231 & {\bf 45} & 187 & 151 & 138 & 231 & {\bf 65} & 203 & 183 & 228 & 273 & {\bf 105} & 315 & 225 & 588 & 273 & {\bf 153} & 379 & 369 \\
$f_{15}$  & 78 & 147 & {\bf 41} & 187 & 135 & 168 & 189 & {\bf 69} & 219 & 159 & 288 & 189 & {\bf 101} & 283 & 191 & 708 & 231 & {\bf 145} & 379 & 277 \\
$f_{16}$  & 18 & 21 & {\bf 9} & 27 & 33 & 18 & 21 & {\bf 9} & 27 & 33 & 18 & 21 & {\bf 9} & 27 & 33 & 18 & 21 & {\bf 9} & 27 & 33 \\
$f_{17}$  & {\bf 79} & 483 & 325 & 619 & 903 & \st{949} & {\bf 777} & 1065 & 1195 & 1491 & 2839 & {\bf 1323} & 1725 & 2123 & 2419 & 6469 & {\bf 1323} & 2077 & 3419 & 2451 \\
$f_{18}$  & 108 & 105 & {\bf 73} & 123 & 145 & 228 & 147 & {\bf 129} & 187 & 209 & 738 & 189 & {\bf 185} & 379 & 395 & 1758 & {\bf 273} & {\bf 273} & 731 & 581 \\
$f_{19}$  & 109 & 231 & {\bf 65} & 155 & 255 & 229 & 231 & {\bf 145} & 475 & 717 & \st{499} & {\bf 231} & 285 & 875 & 1323 & 1369 & {\bf 231} & 449 & 1563 & 1943 \\
$f_{20}$  & 18 & 21 & {\bf 17} & 27 & 33 & 48 & {\bf 21} & 33 & 59 & 33 & {\bf 48} & 63 & 65 & 91 & 95 & 168 & {\bf 63} & 65 & 187 & 219 \\
$f_{21}$  & \st{138} & \st{273} & \st{85} & \st{235} & \st{203} & \st{348} & \st{357} & \st{185} & \st{347} & \st{391} & {\bf 1158} & \st{441} & \st{273} & 1179 & \st{653} & 2748 & \st{525} & {\bf 649} & 1771 & 1839 \\
$f_{22}$  & 228 & {\bf 147} & 241 & 235 & 371 & 888 & 315 & {\bf 305} & 379 & 627 & 2508 & {\bf 315} & 385 & 699 & 627 & 5568 & {\bf 315} & 513 & 1291 & 627 \\
$f_{23}$  & 108 & 273 & {\bf 93} & 299 & 191 & 258 & 399 & {\bf 161} & 411 & 365 & 588 & 441 & {\bf 241} & 699 & 569 & 1608 & 483 & {\bf 401} & 1083 & 957 \\
$f_{24}$  & \st{138} & {\bf 1911} & \st{453} & 3227 & 4515 & \st{1878} & \st{8211} & \st{857} & {\bf 9163} & 11433 & \st{3738} & \st{12285} & \st{1301} & {\bf 15275} & 18503 & \st{5538} & \st{16359} & \st{1745} & {\bf 21147} & 25191 \\
$f_{25}$  & 108 & 567 & {\bf 81} & 379 & 277 & 348 & 819 & {\bf 149} & 859 & 593 & 528 & 819 & {\bf 201} & 1339 & 933 & 678 & 819 & {\bf 269} & 1803 & 1253 \\
    \hline
    \end{tabular}
    \end{scriptsize}\end{center}
    \caption{Results of battery test for $\tau=10^{-3}, 10^{-6}, 10^{-9}$ and $10^{-12}$.
        The columns contain the number of function evaluations required by each integrator
        for each tolerance.
        For each test and tolerance, the best result (least function evaluations) is in bold
        and unsuccessful runs are stricken through.}
    \label{tab:battery}
\end{sidewaystable}

The battery functions were integrated for the relative tolerances
$\tau = 10^{-3}$, $10^{-6}$, $10^{-9}$ and $10^{-12}$ and 
compared to the exact result, computed analytically.
The results are summarized in Table~\ref{tab:battery}, where the
number of required function evaluations for each combination of
integrand, integrator and tolerance are given.
If integration was unsuccessful, the number is stricken through.
If the number of evaluations was the lowest for the given integrand
and tolerance, it is shown in bold face.

\begin{figure}
    \centerline{\epsfig{file=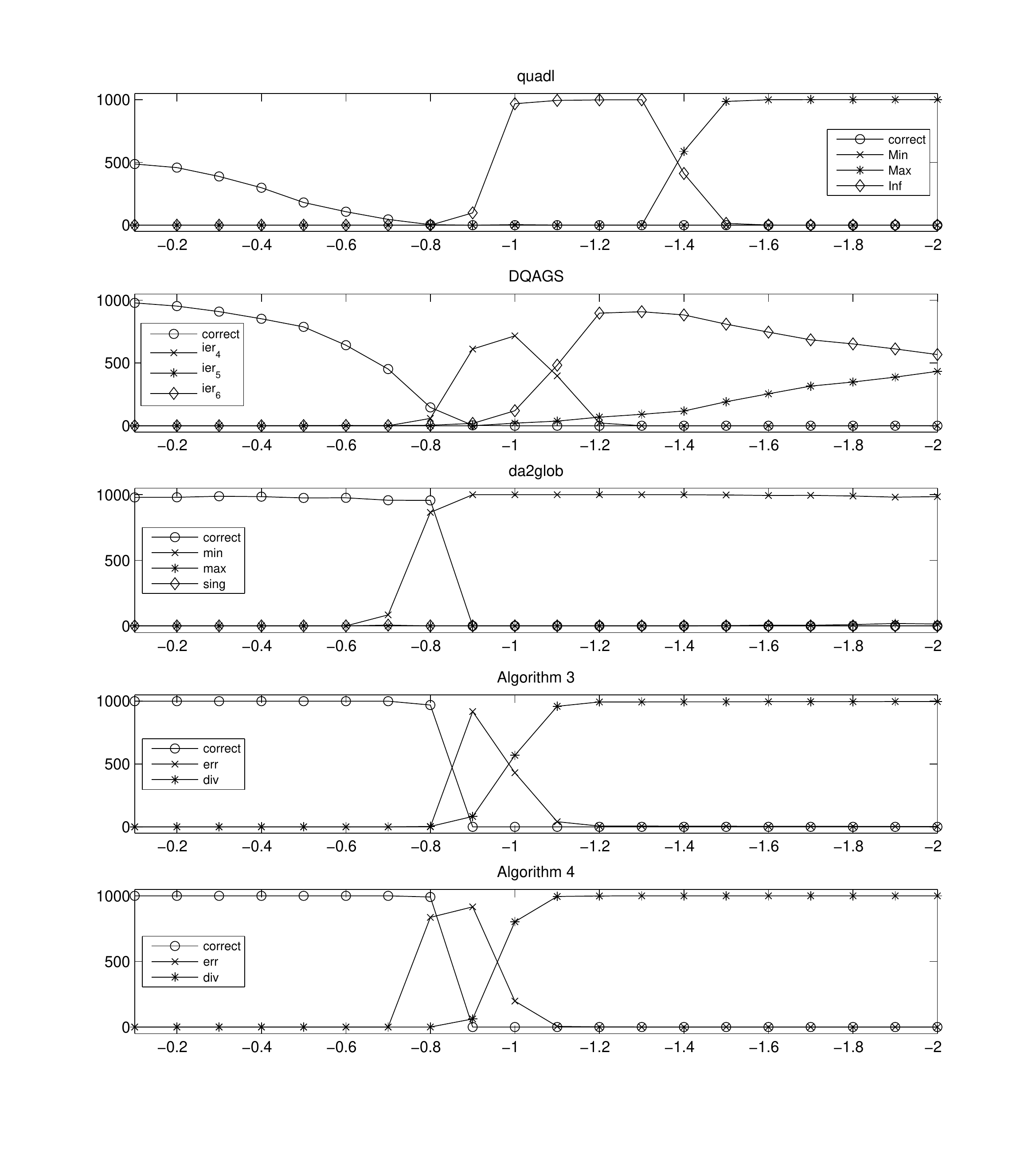,width=0.8\textwidth}}
    \caption{Results of computing $\int_0^1 \left| x - \lambda \right|^\alpha \,\mbox{d}x$
        for $1\,000$ realizations of $\lambda \in [0,1]$ for different 
        values of $\alpha$ ($x$-axis). The curves represent the number of
        correct integrations and the number of times each different warning
        was issued for each value of $\alpha$.}
    \label{fig:divergence}
\end{figure}

Finally, the integrators were tested on the problem
\begin{equation*}
    \int_0^1 \left| x - \lambda \right|^\alpha \,\mbox{d}x, \quad \lambda \in [0,1]
\end{equation*}
for $1\,000$ realizations of the random parameter $\lambda$ and
different values of $\alpha \in [-0.1,-2]$ for a relative\footnote{For $\alpha \leq -1$
an {\em absolute} tolerance of $\tau=10^{-3}$ was used.} tolerance $\tau=10^{-3}$.
Since for $\alpha \leq -1$, the integral diverges and can not be computed
numerically, we are interested in the warnings or errors returned
by the different quadrature routines.
The results are shown in Table~\ref{tab:divergence}
and \fig{divergence}.
For each integrator we give the number of successes and failures
as well as, in brackets, the number of times each possible error
or warning was returned.
The different errors, for each integrator, are:
\begin{itemize}
    \item {\tt quadl}: $(\mathsf{Min}/\mathsf{Max}/\mathsf{Inf})$
        \begin{itemize}
            \item $\mathsf{Min}$: Minimum step size reached; singularity possible.
            \item $\mathsf{Max}$: Maximum function count exceeded; singularity likely.
            \item $\mathsf{Inf}$: Infinite or Not-a-Number function value encountered.
        \end{itemize}
    \item {\tt DQAGS}: $(\mathsf{ier}_1/\mathsf{ier}_2/\mathsf{ier}_3/\mathsf{ier}_4/\mathsf{ier}_5)$
        \begin{itemize}
            \item $\mathsf{ier}_1$: Maximum number of subdivisions allowed
                has been achieved.
            \item $\mathsf{ier}_2$: The occurrence of roundoff error was detected,
                preventing the requested tolerance from being achieved.
                The error may be under-estimated.
            \item $\mathsf{ier}_3$: Extremely bad integrand behavior somewhere
                in the interval.
            \item $\mathsf{ier}_4$: The algorithm won't converge due to roundoff
                error detected in the extrapolation table. It is presumed that
                the requested tolerance cannot be achieved, and that the
                returned result is the best which can be obtained.
            \item $\mathsf{ier}_5$: The integral is probably divergent
                or slowly convergent.
        \end{itemize}
    \item {\tt da2glob}: $(\mathsf{noise}/\mathsf{min}/\mathsf{max}/\mathsf{sing})$
        \begin{itemize}
            \item $\mathsf{noise}$: The requested tolerance is below the noise level
                of the problem. Required tolerance may not be met.
            \item $\mathsf{min}$: Interval too small.
                Required tolerance may not be met.
            \item $\mathsf{max}$: Maximum number of function evaluations.
                Required tolerance may not be met.
            \item $\mathsf{sing}$: Singularity probably detected.
                Required tolerance may not be met.
        \end{itemize}
    \item Algorithms~\ref{alg:final_naive} and \ref{alg:final_refined}:
        $(\mathsf{err}/\mathsf{div})$
        \begin{itemize}
            \item $\mathsf{err}$: The final error estimate is larger than the
                required tolerance.
            \item $\mathsf{div}$: The integral is divergent.
        \end{itemize}
\end{itemize}
Thus, the results for {\tt DQAGS} at $\alpha = -0.8$ should be read
as the algorithm returning 146 correct and 854 false (requested
tolerance not satisfied) results and having returned the error $\mathsf{ier}_3$
(bad integrand behavior) 58 times and the error $\mathsf{ier}_5$
(probably divergent integral) 4 times.

\begin{sidewaystable}
    \small
    \begin{center}\begin{scriptsize}
        \begin{tabular}{l|c|c|c|c|c}
            $\alpha$ & {\tt quadl} & {\tt DQAGS} & {\tt da2glob} & Algorithm~\ref{alg:final_refined} & Algorithm~\ref{alg:final_naive} \\ \hline
$\alpha = -0.1$ & 487 / 513 (0/0/0) & 979 / 21 (0/0/0/0/0) & 979 / 21 (0/0/0/0) & 1000 / 0 (0/0) & 1000 / 0 (0/0)  \\ \hline
$\alpha = -0.2$ & 459 / 541 (0/0/0) & 953 / 47 (0/0/0/0/0) & 980 / 20 (0/0/0/0) & 1000 / 0 (0/0) & 1000 / 0 (0/0)  \\ \hline
$\alpha = -0.3$ & 388 / 612 (0/0/0) & 909 / 91 (0/0/0/0/0) & 988 / 12 (0/0/0/0) & 1000 / 0 (0/0) & 1000 / 0 (0/0)  \\ \hline
$\alpha = -0.4$ & 298 / 702 (0/0/0) & 852 / 148 (0/0/0/0/0) & 985 / 15 (0/0/0/0) & 1000 / 0 (0/0) & 1000 / 0 (0/0)  \\ \hline
$\alpha = -0.5$ & 181 / 819 (0/0/0) & 788 / 212 (0/0/0/0/1) & 975 / 25 (0/0/0/0) & 1000 / 0 (0/0) & 1000 / 0 (0/0)  \\ \hline
$\alpha = -0.6$ & 107 / 893 (0/0/0) & 642 / 358 (0/0/0/0/2) & 976 / 24 (0/0/0/0) & 1000 / 0 (0/0) & 1000 / 0 (0/0)  \\ \hline
$\alpha = -0.7$ & 46 / 954 (0/0/0) & 451 / 549 (0/0/0/0/1) & 958 / 42 (0/84/0/6) & 1000 / 0 (0/0) & 1000 / 0 (0/0)  \\ \hline
$\alpha = -0.8$ & 4 / 996 (0/0/0) & 146 / 854 (0/0/58/0/4) & 957 / 43 (0/865/0/0) & 998 / 2 (836/0) & 970 / 30 (0/3)  \\ \hline
$\alpha = -0.9$ & 0 / 1000 (0/0/98) & 3 / 997 (0/0/611/0/18) & 0 / 1000 (0/1000/0/0) & 0 / 1000 (916/62) & 0 / 1000 (918/82)  \\ \hline
\hline
$\alpha = -1.0$ & 0 / 1000 (4/0/968) & 0 / 1000 (0/0/716/20/119) & 0 / 1000 (0/1000/0/0) & 0 / 1000 (198/802) & 0 / 1000 (431/569)  \\ \hline
$\alpha = -1.1$ & 0 / 1000 (1/0/994) & 0 / 1000 (0/0/397/36/482) & 0 / 1000 (0/1000/0/0) & 0 / 1000 (5/995) & 0 / 1000 (40/960)  \\ \hline
$\alpha = -1.2$ & 0 / 1000 (1/0/998) & 0 / 1000 (0/0/21/68/897) & 0 / 1000 (0/1000/0/0) & 0 / 1000 (1/999) & 0 / 1000 (6/994)  \\ \hline
$\alpha = -1.3$ & 0 / 1000 (1/0/999) & 0 / 1000 (0/0/0/91/908) & 0 / 1000 (0/1000/0/0) & 0 / 1000 (0/1000) & 0 / 1000 (6/994)  \\ \hline
$\alpha = -1.4$ & 0 / 1000 (0/588/412) & 0 / 1000 (0/0/0/117/883) & 0 / 1000 (0/1000/0/0) & 0 / 1000 (0/1000) & 0 / 1000 (5/995)  \\ \hline
$\alpha = -1.5$ & 0 / 1000 (0/986/14) & 0 / 1000 (0/0/0/191/809) & 0 / 1000 (0/998/2/0) & 0 / 1000 (0/1000) & 0 / 1000 (5/995)  \\ \hline
$\alpha = -1.6$ & 0 / 1000 (0/999/1) & 0 / 1000 (0/0/0/254/746) & 0 / 1000 (0/994/6/0) & 0 / 1000 (0/1000) & 0 / 1000 (4/996)  \\ \hline
$\alpha = -1.7$ & 0 / 1000 (0/1000/0) & 0 / 1000 (0/0/0/315/685) & 0 / 1000 (0/995/5/0) & 0 / 1000 (0/1000) & 0 / 1000 (4/996)  \\ \hline
$\alpha = -1.8$ & 0 / 1000 (0/1000/0) & 0 / 1000 (0/0/0/348/652) & 0 / 1000 (0/990/10/0) & 0 / 1000 (0/1000) & 0 / 1000 (4/996)  \\ \hline
$\alpha = -1.9$ & 0 / 1000 (0/1000/0) & 0 / 1000 (0/0/1/387/612) & 0 / 1000 (0/981/19/0) & 0 / 1000 (0/1000) & 0 / 1000 (3/997)  \\ \hline
$\alpha = -2.0$ & 0 / 1000 (0/1000/0) & 0 / 1000 (0/0/0/433/567) & 0 / 1000 (0/985/15/0) & 0 / 1000 (0/1000) & 0 / 1000 (3/997)  \\ \hline
        \end{tabular}
    \end{scriptsize}\end{center}
    \caption{Results of computing $\int_0^1 \left| x - \lambda \right|^\alpha \,\mbox{d}x$
        for $1\,000$ realizations of $\lambda \in [0,1]$ for different $\alpha$.
        The columns contain the number of correct/incorrect integrations as well
        as the number of times (in brackets) the different errors or warnings
        of each algorithm were returned. Thus, the results for {\tt DQAGS} at $\alpha = -0.8$ should be read
        as the algorithm returning 146 correct and 854 false 
        results and having returned the error $\mathsf{ier}_3$
        58 times and the error $\mathsf{ier}_5$
        4 times.}
    \label{tab:divergence}
\end{sidewaystable}

\section{Discussion}
\label{sec:discussion}

As can be seen from the results in Table~\ref{tab:LKtest}, for
the integrands in the Lyness-Kaganove test, both new
algorithms presented in Section~\ref{sec:algorithm} are clearly
more reliable than {\tt quadl} and {\tt DQAGS}.
{\small MATLAB}'s {\tt quadl} performs best for high precision requirements
(small tolerances, best results for $\tau = 10^{-9}$),
yet still fails often without warning.
{\small QUADPACK}'s {\tt DQAGS} does better for low precision requirements
(large tolerances, best results for $\tau = 10^{-3}$), yet also fails often, 
more often than not with a warning.

\begin{figure}
    \begin{center}
        \epsfig{file=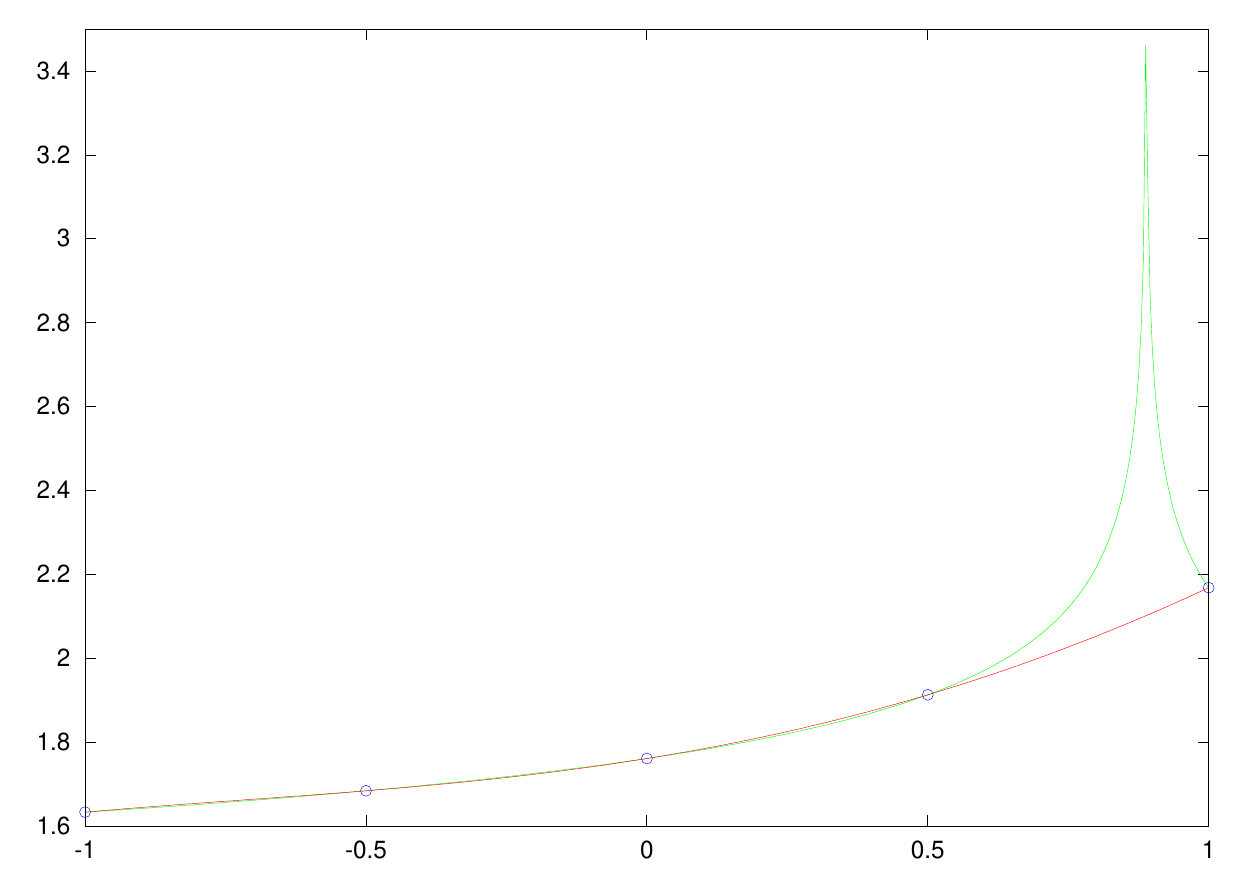,width=0.45\textwidth}
        \epsfig{file=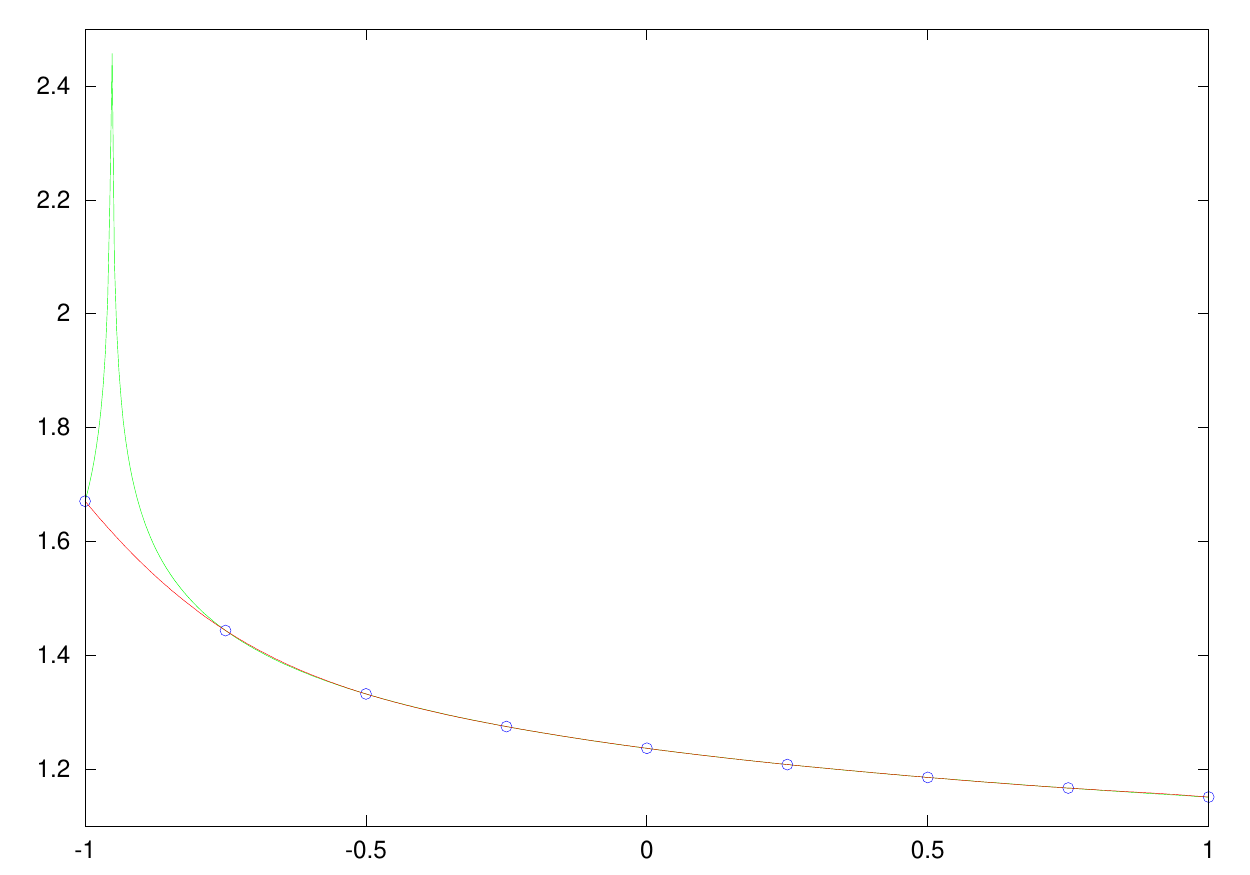,width=0.45\textwidth}
        \epsfig{file=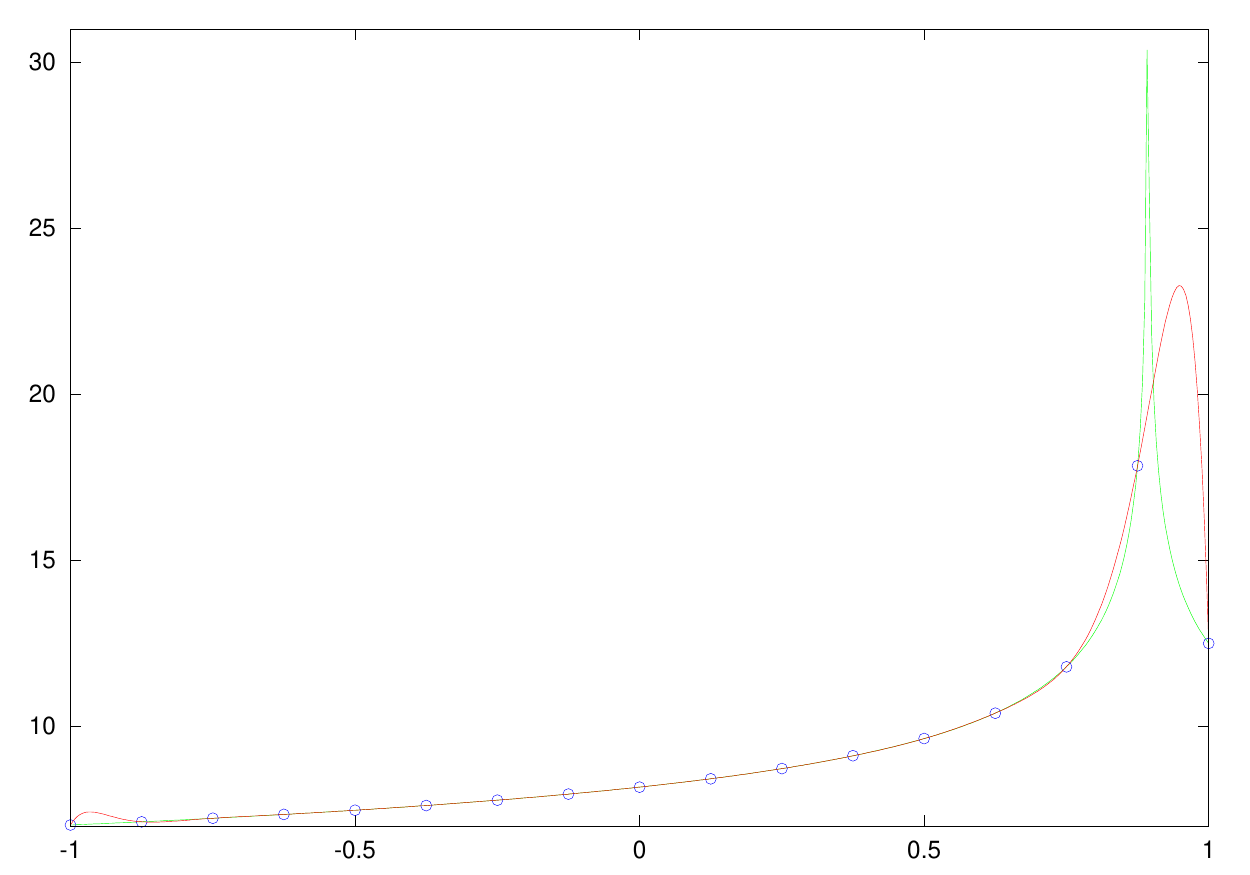,width=0.45\textwidth}
    \end{center}
    \caption{Some cases in which the 5, 9 and 17-point rules used in
        {\tt da2glob} fail to detect a singularity (as in \eqn{val_sing}). The
        assumed integrand (red) is sufficiently smooth such that
        its higher-degree coefficients, from which the error estimate
        is computed, are near zero.}
    \label{fig:nonsmooth}
\end{figure}

Espelid's {\tt da2glob} does significantly better, with only
a few failures at $\tau = 10^{-3}$ (without warning) and
a large number of failures for \eqn{val_oscill} at
$\tau = 10^{-12}$, albeit all of them with prior warning.
The former are due to the error estimate not detecting specific
features of the integrand due to the interpolating polynomial
looking smooth, when it is, in fact, singular (see Figure~\ref{fig:nonsmooth}).
The latter were due to the integral estimates being affected by noise
(most often in the 17-point rule),
which was, in all cases, detected and warned against by the algorithm.

The new algorithms fail only for \eqn{val_sing} at small tolerances
since the integral becomes numerically impossible to evaluate
(there are not sufficient machine numbers near the singularity
to properly resolve the integrand),
for which a warning is issued.
This problem is shared by the other integrators as well.
Algorithm~\ref{alg:final_naive} also fails a few times
when integrating \eqn{val_4peak}.
In all such cases, one of the peaks was missed completely
by the lower-degree rules, giving the appearance of a flat curve.
Both algorithms also failed a few times on \eqn{val_oscill}
in cases where the resulting integral was several orders of magnitude smaller than
the function values themselves, making the required tolerance practically
un-attainable.

Whereas the new algorithms out-perform the others in terms of reliability,
they do so at a cost of a higher number of function evaluations.
On average, Algorithm~\ref{alg:final_naive} uses about {\em twice as many}
function evaluations as {\tt da2glob}, whereas Algorithm~\ref{alg:final_refined}
uses roughly {\em six times} as many.

\begin{figure}
    \begin{center}
        \epsfig{file=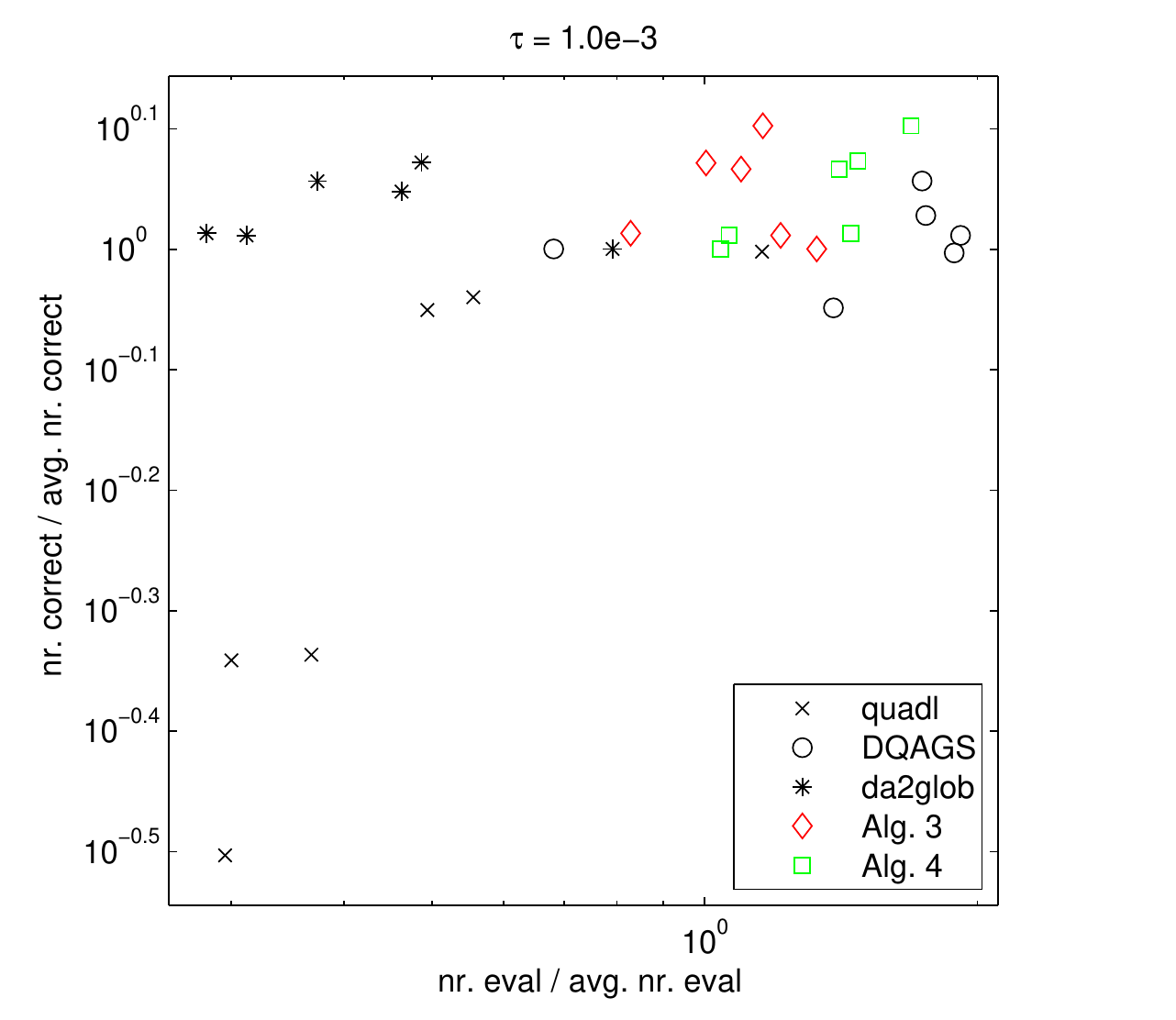,width=0.48\textwidth}
        \epsfig{file=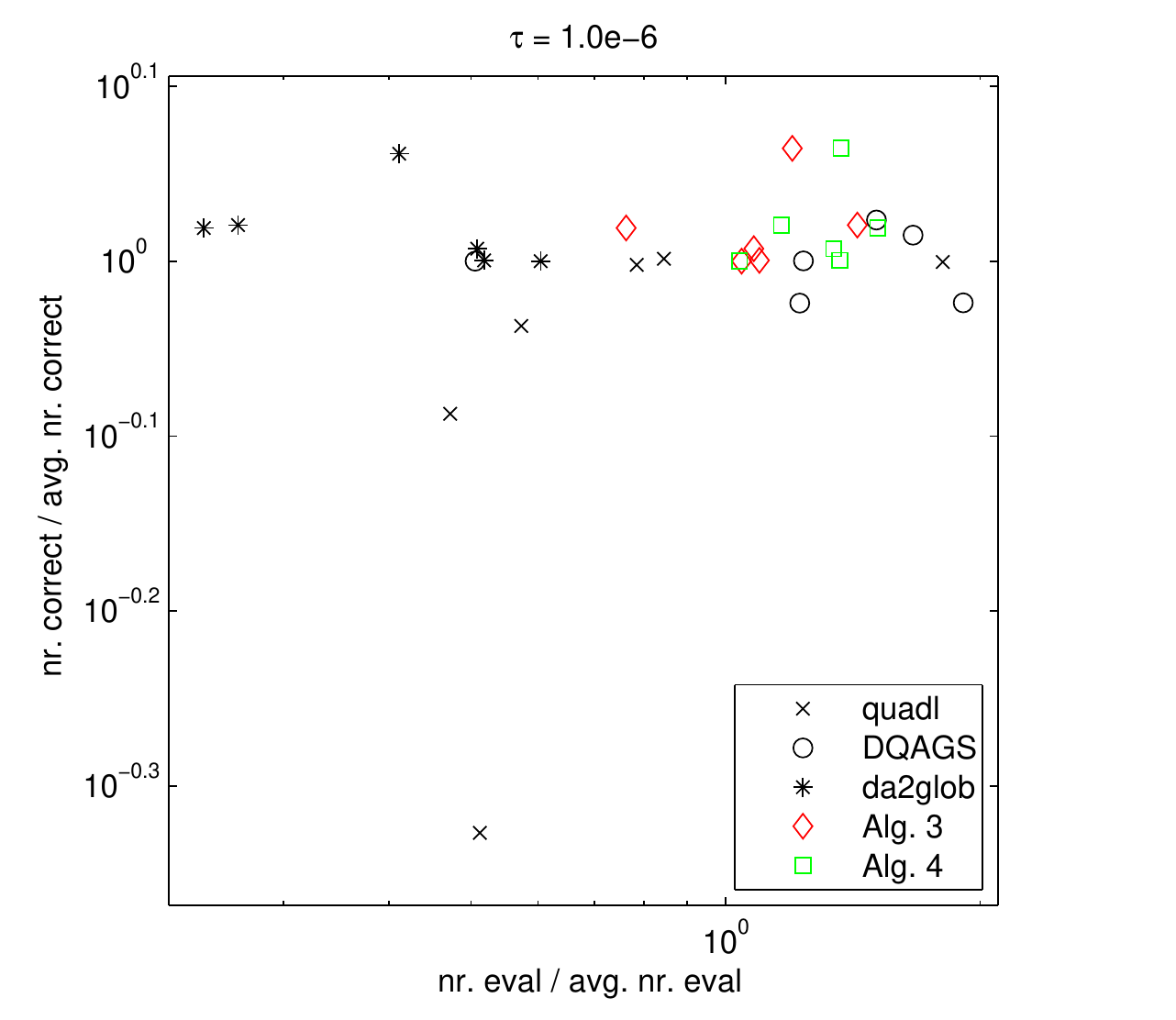,width=0.48\textwidth}
        \epsfig{file=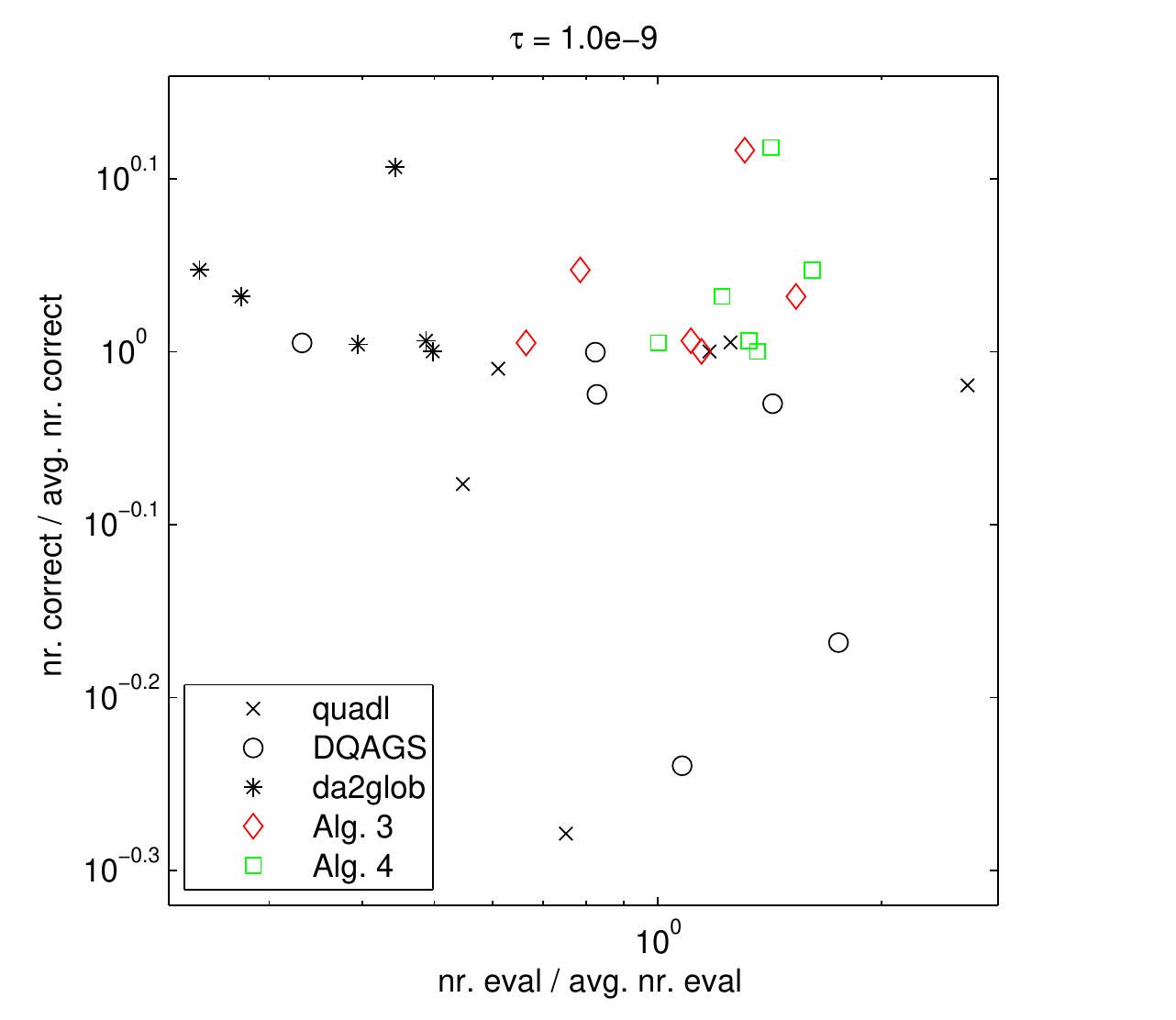,width=0.48\textwidth}
        \epsfig{file=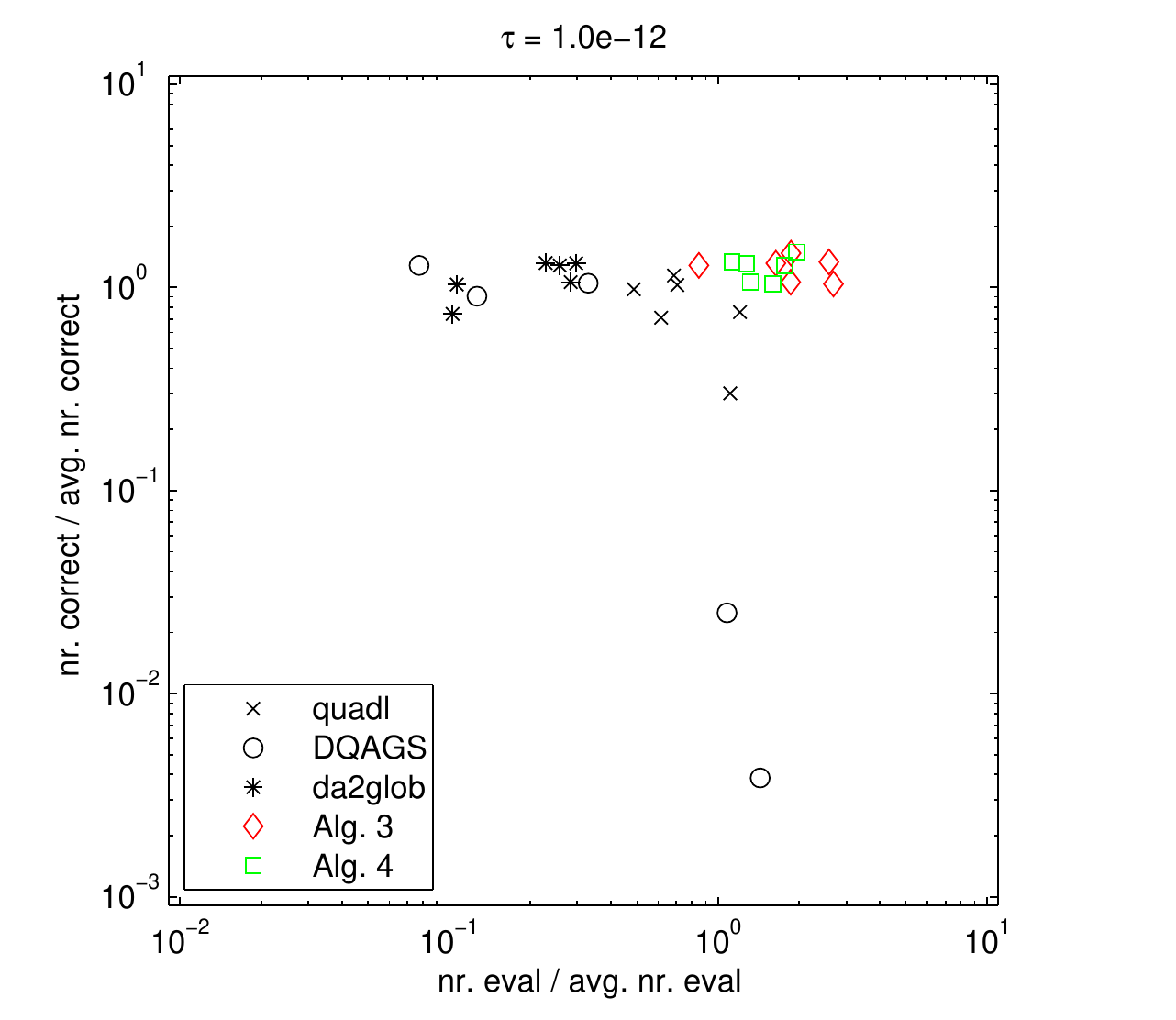,width=0.48\textwidth}
    \end{center}
    \caption{Scatter-plots of the results of the Lyness-Kaganove
        test-suite for each tolerance. Each point represents one of
        the test functions (Equations~\ref{eqn:val_sing} to \ref{eqn:val_oscill}).
        Its location is determined by the relative number of function evaluations
        (on the $x$-axis) and the relative number of correct evaluations
        (on the $y$-axis).}
    \label{fig:results}
\end{figure}

The results are summarized in Figure~\ref{fig:results}.
The plots, for each tolerance and integrand, show where each
algorithm stand in terms of relative reliability and relative efficiency.
If we divide the plots into four regions

\centerline{\epsfig{file=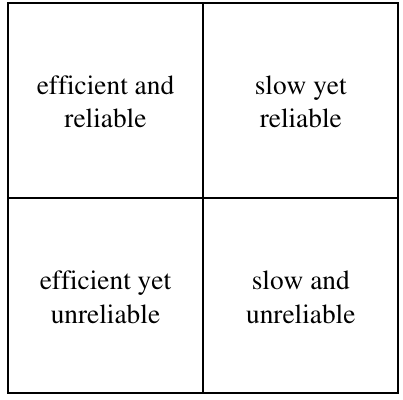,width=0.2\textwidth}}

\noindent we can see that whereas {\tt da2glob} is clearly efficient and reliable,
the two new algorithms are slightly more reliable, yet slow.
The results for {\tt quadl} and {\small QUADPACK}'s {\tt DQAGS}
are scattered over all four regions.

This trend is also visible in the results of the ``battery'' test
(Table~\ref{tab:battery}).
Algorithms~\ref{alg:final_naive} and \ref{alg:final_refined}
fail on $f_{21}$ for all but the highest and second-highest tolerances
respectively, since the third peak at $x=0.6$ is missed completely.
{\tt da2glob} also does quite well, failing on $f_{21}$ at the same tolerances
and for the same reasons as the new algorithms and on $f_{24}$
for $\tau < 10^{-3}$, using, however, in almost all cases,
less function evaluations than Algorithms~\ref{alg:final_naive}
or \ref{alg:final_refined}.

\begin{figure}
    \begin{center}
        a~\epsfig{file=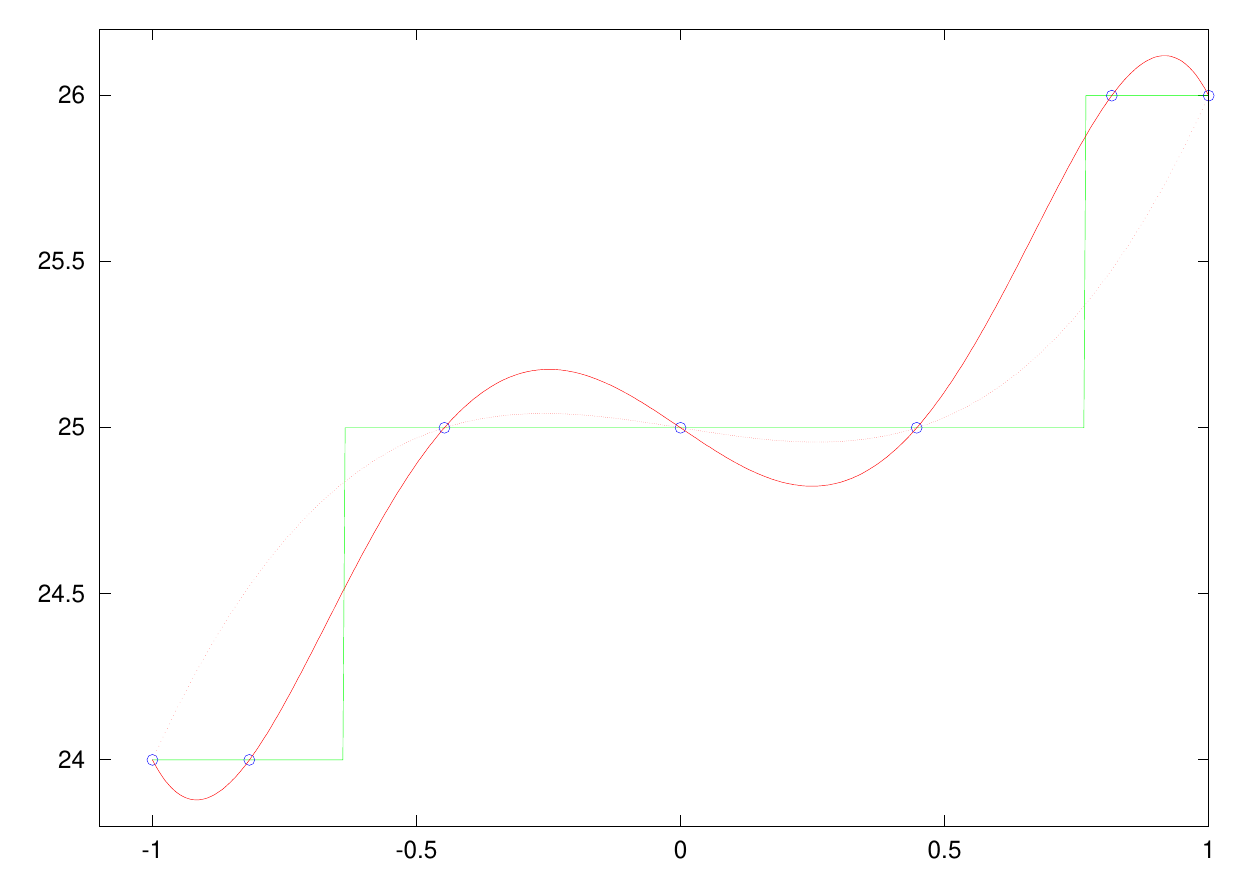,width=0.45\textwidth}
        b~\epsfig{file=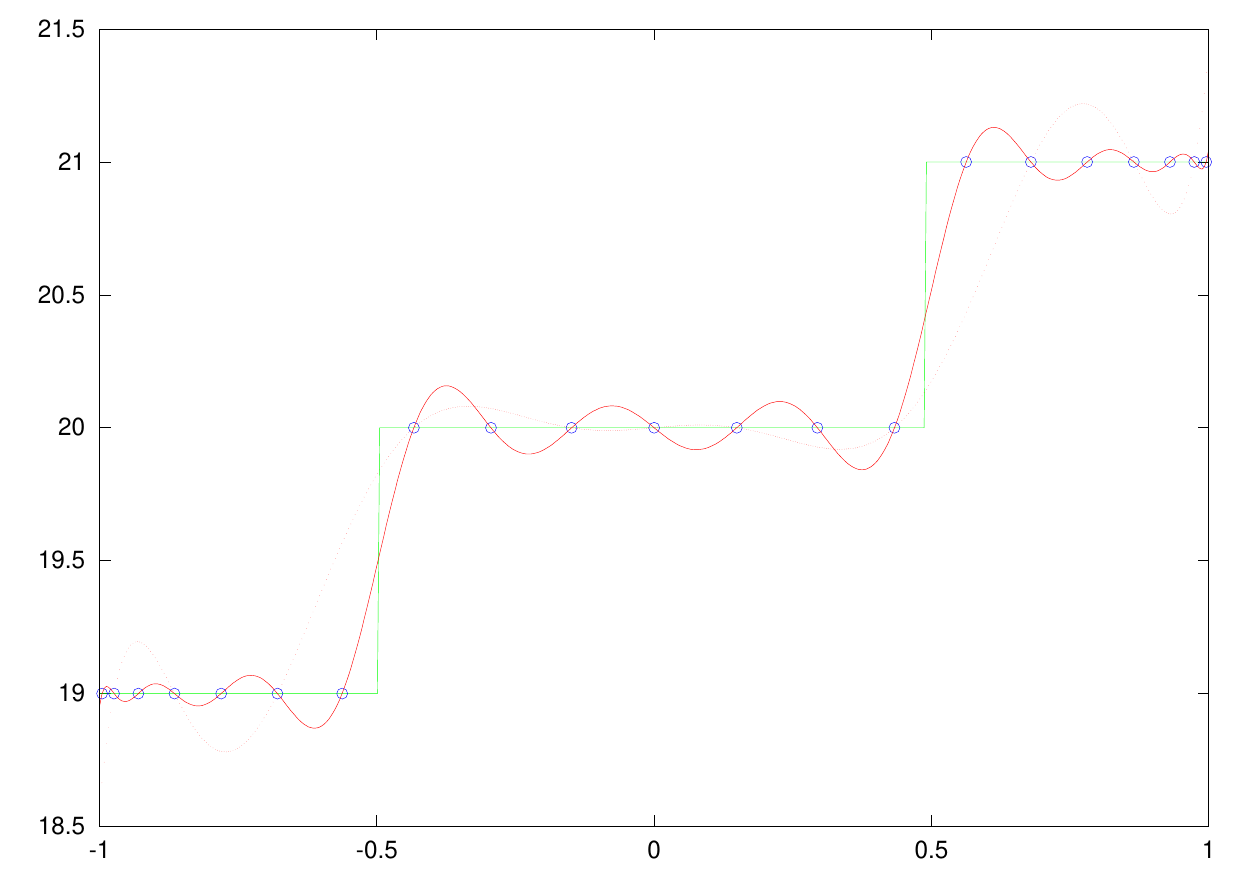,width=0.45\textwidth}
        c~\epsfig{file=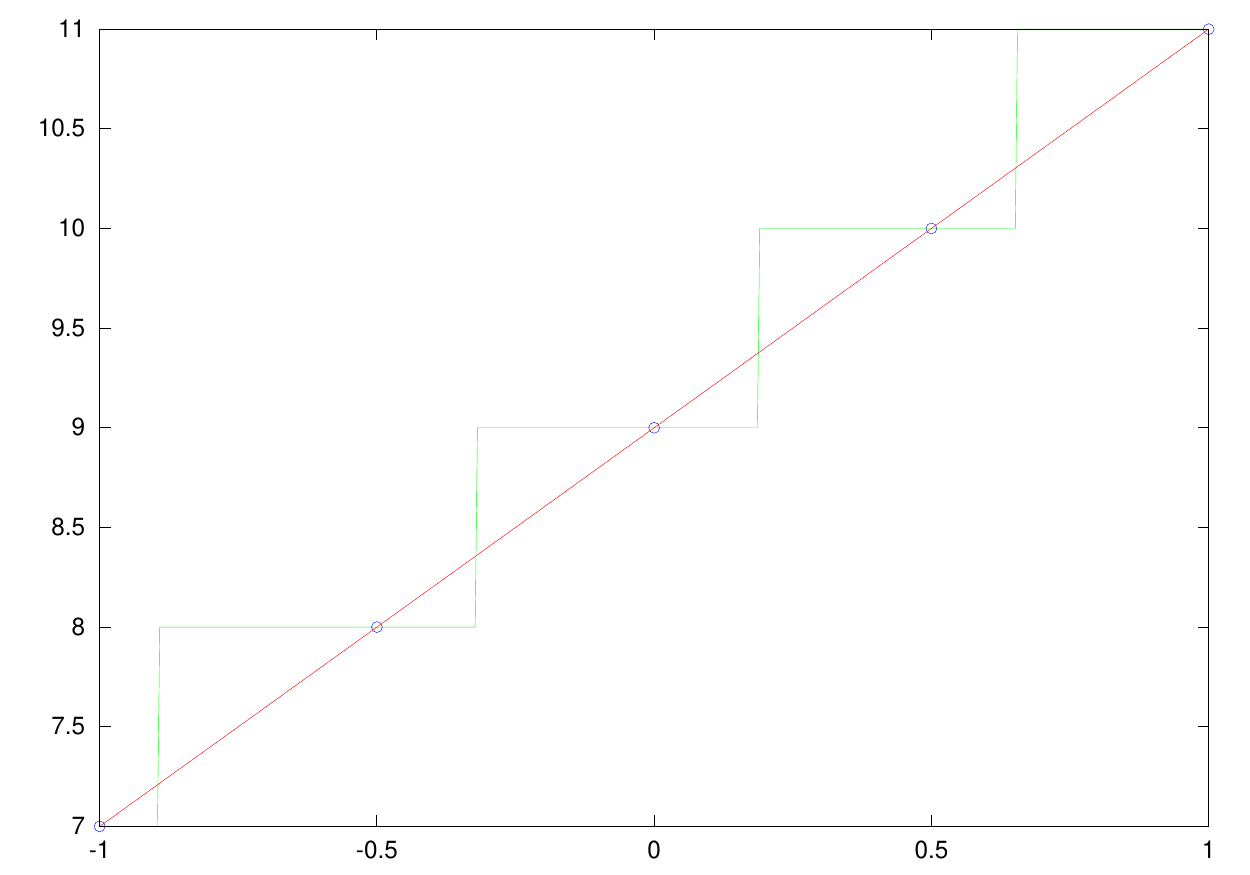,width=0.45\textwidth}
    \end{center}
    \caption{Intervals in which the error estimates of {\tt quadl} (a),
        {\tt DQAGS} (b) and {\tt da2glob} (c) failed for $f_{24}$.
        The interpolatory polynomials used to compute the error
        estimates are shown in red.}
    \label{fig:waldvogel}
\end{figure}

It is interesting to note that {\tt quadl}, {\tt DAQGS} and {\tt da2glob}
all failed to integrate $f_{24}$ for tolerances $\tau<10^{-3}$.
A closer look at the intervals that caused each algorithm to fail
(see Figure~\ref{fig:waldvogel}) reveals that in all cases, multiple
discontinuities in the same interval caused the error estimate to be
accidentally small, leading the algorithms to erroneously assume convergence.
This is not a particularity of the interval chosen: if we integrate
\begin{equation}
    \label{eqn:waldvogel}
    \int_0^\lambda \lfloor e^x \rfloor \,\mbox{d}x, \quad \lambda \in [2.5,3.5]
\end{equation}
for $1\,000$ realizations of the random parameter $\lambda$ for
$\tau = 10^{-6}$ using these three integrators, they fail on
$894$, $945$ and $816$ of the cases respectively.
Both Algorithms~\ref{alg:final_naive} and \ref{alg:final_refined}
succeed in all cases.

One could argue that the different error estimates are all equally
good and the different ratios of reliability vs.\ 
efficiency are only due to their parameterization, \ie their scaling
of the error estimate.
We can test this hypothesis by scaling the error estimates of all
five algorithms by the smallest value\footnote{For simplicity, we consider
only values in units of the next-closest power of 10.} $\rho$ such that
the 1\,000 runs at $\tau=10^{-3}$ for \eqn{val_sing}
produce no incorrect results:

\begin{center}\begin{small}\begin{tabular*}{0.9\textwidth}{@{\extracolsep{\fill}}l|c|c|c|c|c|c|c|c|c|c}
    $\tau=10^{-3}$ & \multicolumn{2}{c|}{\tt quadl} & \multicolumn{2}{c|}{\tt DQAGS} & \multicolumn{2}{c|}{\tt da2glob} & \multicolumn{2}{c|}{Algorithm 4} & \multicolumn{2}{c}{Algorithm 3} \\
    $f(x)$ & $\rho$ & $n_\mathsf{eval}$ & $\rho$ & $n_\mathsf{eval}$ & $\rho$ & $n_\mathsf{eval}$ & $\rho$ & $n_\mathsf{eval}$ & $\rho$ & $n_\mathsf{eval}$ \\ \hline
    Eqn (\ref{eqn:val_sing}) & 40\,000 & 783.18 & --- & --- & 50 & 229.40 & 0.006 & 319.24 & 0.02 & 243.85 \\ \hline
\end{tabular*}\end{small}\end{center}

\noindent For {\tt quadl}, a scaling of $\rho=40\,000$ was necessary,
resulting in an 8-fold increase in the number of required evaluations and
for {\small QUADPACK}'s {\tt DQAGS}, no amount of scaling produced
the correct result.
For {\tt da2glob}, Algorithm~\ref{alg:final_naive} and 
Algorithm~\ref{alg:final_refined},
only moderate scaling
was necessary, bringing the number of required function evaluations
much closer to each other.
There is, therefore, for {\tt da2glob} a potential for tuning empirical scaling
factors towards more reliability, yet at the cost of close to
the same number of function evaluations as Algorithms~\ref{alg:final_naive}
and \ref{alg:final_refined}, at least for \eqn{val_sing}.

This approach, however, breaks down completely for \eqn{waldvogel}.
As can be seen in \fig{waldvogel}, the error estimates which cause
the individual algorithms to fail are not merely small: they are zero,
and hence no amount of scaling will fix them.
The two new algorithms are therefore not only more reliable due to
a stricter (or more pessimistic) scaling of the error estimate,
but mainly due to the fundamentally different type of error estimate,
which is less prone to accidentally small estimates (see \cite{ref:Gonnet2009}).

In the final test evaluating divergent integrals (Table~\ref{tab:divergence}),
{\tt quadl} fails to 
distinguish between divergent and non-divergent integrals, reporting
that a non-numerical value was encountered for $-1.0 \geq \alpha \geq -1.4$
and then aborting after the maximum $10\,000$ function evaluations\footnote{
This termination criteria had been disabled for the previous tests.}
for $\alpha < 1.4$. 
For $\alpha < -1.0$, {\tt DQAGS} reports the integral to be either
subject to too much rounding error or divergent.
The latter correct result was returned in more than half of the cases tested.
In most cases where $\alpha < -0.8$, {\tt da2glob} aborted, reporting
that the smallest interval size had been reached, or, in some cases,
that the maximum number of evaluations (by default $10\,000$) had been
exceeded.
All these cases were accompanied by an additional warning that a singularity
had probably been detected.
For $\alpha = -0.8$, Algorithm~\ref{alg:final_refined} fails with a warning
that the required tolerance was not met and as of $\alpha < -1.1$
both Algorithms~\ref{alg:final_refined} and \ref{alg:final_naive}
abort, correctly, after deciding that the integral is divergent.

We conclude that the new Algorithms~\ref{alg:final_refined} and
\ref{alg:final_naive}, presented herein, are {\em more reliable}
than {\small MATLAB}'s {\tt quadl}, {\small QUADPACK}'s {\tt DQAGS}
and Espelid's {\tt da2glob}.
This higher reliability is not due to a stricter scaling of the error, but
to a new type of error estimator which avoids most of the problems
observed in the other algorithms.
This increased reliability comes at a cost of two to six times higher
number of function evaluations required for complicated integrands
such as those in Equations~\ref{eqn:val_sing} to \ref{eqn:val_oscill}.

The tradeoff between reliability and efficiency should, however,
be of little concern
in the context of automatic or general-purpose quadrature routines,
which should work reliably for {\em any} type of integrand.
Most modifications which increase efficiency usually rely on making
certain assumptions on the integrand, \eg smoothness, continuity, 
non-singularity, monotonically decaying coefficients, etc\dots \xspace
If, however, the user knows enough about his or her integrand as to 
know that these assumptions are indeed valid and therefore that the
algorithm will not fail, then he or she knows enough about the integrand
as to {\em not} have to use a general-purpose quadrature routine and,
if efficiency is crucial, should consider integrating it by a 
specialized routine or even trying to do so analytically.

In making any assumptions for the user, we would be making two mistakes:
\begin{enumerate}
    \item The increase in efficiency would {\em reward} users who,
        despite knowing enough about their integrand to trust
        the quadrature rule, have not made the effort
        to look for a specialized or less general routine,
    \item The decrease in reliability {\em punishes} users who have
        turned to a general-purpose quadrature routine because they
        knowingly {\em can not} make any assumptions regarding their
        integrand.
\end{enumerate}
It is for this reason that we should have no qualms whatsoever in sacrificing a bit
of efficiency on some special integrands for much more reliability on
tricky integrands for which we know, and can therefore assume, nothing.

Ideally, software packages such as Matlab or libraries such as the
Gnu Scientific Library \cite{ref:Galassi2009} should offer both 
heavy-duty quadrature routines such as Algorithms~\ref{alg:final_naive}
and \ref{alg:final_refined} presented herein, alongside other efficient
routines, such as {\tt da2glob} \cite{ref:Espelid2007} for less tricky
integrands, as well as even more efficient methods (\eg the exponentially
convergent integrator proposed by \citeN{ref:Waldvogel2009}) for
analytic (continuous and smooth) integrands, allowing the user to choose
according to his or her specific needs or level of understanding
regarding the integrand.

\section{Acknowledgements}
I would like to thank Walter Gander who supervised my PhD thesis which
resulted in this work.
Further thanks go to J\"org Waldvogel, Fran\c{c}ois Cellier,
Gradimir Milovanovi\'c, Aleksandar Cvetkovi\'c,
Marija Stani\'c, Geno Nikolov and Borislav Bojanov who,
through their collaboration in the Swiss National Science Foundation (SNSF)
SCOPES Project (Nr.~IB7320-111079/1, 2005--2008) ``New Methods for Quadrature'',
provided not only valuable answers, but also the right questions,
as well as to Gaston Gonnet for his input on several topics regarding
the algorithms and the manuscript itself.
Finally, I would like to thank the anonymous reviewers who's
corrections and suggestions, especially those regarding the 
graphical representation of the results, made this a better
manuscript.

\bibliography{quad}

\begin{received}
...
\end{received}
\end{document}